\documentclass{amsart}
\usepackage{amssymb, amsmath, amscd, latexsym, mathrsfs, eufrak, lscape, tabls}
\input xy
\xyoption{all}
\newtheorem{thm}{Theorem}[section]
\newtheorem{cor}[thm]{Corollary}
\newtheorem{lem}[thm]{Lemma}

\newtheorem{prop}[thm]{Proposition}
\newtheorem{defn}[thm]{Definition}

\newtheorem{rem}[thm]{Remark}

\newtheorem{claim}[thm]{Claim}
\newtheorem{ex}{Example}[section]
\numberwithin{equation}{section}

\newtheorem{bst}{Bisymmetric Triple}[section]
\newtheorem{symp}{Symmetric Pair}[section]
\newenvironment{pr}{\textit{Proof: }}

\newcommand{\al}{\alpha}
\newcommand{\be}{\beta}
\newcommand{\la}{\lambda}

\newcommand{\ga}{\gamma}

\newcommand{\de}{\delta}

\newcommand{\li}{\medskip}

\newcommand{\reals}{\mathbb{R}}

\newcommand{\complex}{\mathbb{C}}

\newcommand{\rationals}{\mathbb{Q}}

\newcommand{\mfg}{\mathfrak{g}}

\newcommand{\mfk}{\mathfrak{k}}
\newcommand{\mfl}{\mathfrak{l}}
\newcommand{\mfm}{\mathfrak{m}}
\newcommand{\mfn}{\mathfrak{n}}

\newcommand{\mfp}{\mathfrak{p}}

\def\bar{\begin{array}}
\def\ear{\end{array}}
\def\sbar{\begin{subarray}}
\def\sear{\end{subarray}}
\def\beq{\begin{equation}}
\def\eeq{\end{equation} }
\def\beqar{\begin{eqnarray}}
\def\eeqar{\end{eqnarray}}
\def\bal{\begin{align}}
\def\eal{\end{align}}
\def\bfig{\begin{figure}}
\def\efig{\end{figure}}
\def\bc{\begin{center}}
\def\ec{\end{center}}
\def\btab{\begin{table}}
\def\etab{\end{table}}
\def\bland{\begin{landscape}}
\def\eland{\end{landscape}}

\def\bproof{\begin{pr}}
\def\eproof{\end{pr}}
\def\bprop{\begin{prop}}
\def\eprop{\end{prop}}
\def\bthm{\begin{thm}}
\def\ethm{\end{thm}}
\def\blem{\begin{lem}}
\def\elem{\end{lem}}
\def\brem{\begin{rem}}
\def\erem{\end{rem}}
\def\bcor{\begin{cor}}
\def\ecor{\end{cor}}
\def\bex{\begin{ex}}
\def\eex{\end{ex}}
\def\bdfn{\begin{defn}}
\def\edfn{\end{defn}}
\def\bclaim{\begin{claim}}
\def\eclaim{\end{claim}}
\def\bbst{\begin{bst}}
\def\ebst{\end{bst}}
\def\bsymp{\begin{symp}}
\def\esymp{\end{symp}}

\newlength{\myVSpace}
\newcommand\xstrut{\raisebox{-.9\myVSpace}{\rule{0pt}{\myVSpace}}}% to change the height of rows in arrays

\makeatletter
\def\table{\@ifnextchar[{\table@i}{\table@i[\fps@table]}}
\def\table@i[#1]{\@float{table}[#1]\footnotesize}
\makeatother

\begin{document}

\title[Some Einstein Homogeneous Riemannian Fibrations]
{Some Einstein Homogeneous Riemannian Fibrations}

\author{F\'{a}tima Ara\'{u}jo}

\address{School of Mathematics, The University of Edinburgh, JCMB, The Kings Buildings, Edinburgh, EH9 3JZ}

\email{m.d.f.araujo@sms.ed.ac.uk}

%\thanks{}
%\subjclass{}%
%\keywords{}%

\date{}%
%\dedicatory{}%
%\commby{}%
% ----------------------------------------------------------------
%\begin{abstract}
%\end{abstract}
\maketitle
% ----------------------------------------------------------------

\begin{abstract}
We study the existence of projectable $G$-invariant Einstein metrics on the total space of $G$-equivariant fibrations $M=G/L\rightarrow G/K$, for a compact connected semisimple Lie group $G$. We obtain necessary conditions for the existence of such Einstein metrics in terms of appropriate Casimir operators, which is a generalization of the result by Wang and Ziller about Einstein normal metrics. We describe binormal Einstein metrics which are the orthogonal sum of the normal metrics on the fiber and on the base. The special case when the restriction to the fiber and the projection to the base are also Einstein is also considered. As an application, we prove the existence of a non-standard Einstein invariant metric on the Kowalski $n$-symmetric spaces.
\end{abstract}

\li

\let\thefootnote\relax\footnotetext{\emph{MSC-class:} 53C25; 53C30; 53C20; 53C35.}

\let\thefootnote\relax\footnotetext{\emph{Keywords and phrases:} Einstein metrics, homogeneous space, fibration, totally geodesic.}
   
\let\thefootnote\relax\footnotetext{\emph{Acknowledgements:} I would like to thank Dmitri Alekseevsky for his useful advice and enlightening discussions. The author was sponsored by Funda\c{c}\~{a}o para a Ci\^{e}ncia e a Tecnologia SFRH/BD/12267/2003.}

\maketitle

\section{Introduction}\label{intro}

We describe a class of invariant Einstein metrics on a homogeneous manifold $M=G/L$ of a compact connected semisimple Lie group $G$, which is consistent with a homogeneous fibration $G/L\rightarrow G/K$.

A Riemannian metric $g$ is said to be Einstein if its Ricci curvature satisfies the Einstein equation $Ric=Eg$, for some constant $E$. The Einstein equation is a system of partial differential equations of second order, which is in general unmanageable. A few results about Einstein metrics are known in the general case and many results are known under some extra assumptions (special type of metrics or metric with large isometry group). For example, there are deep results about K\"{a}hler-Einstein (\cite{Yau}, \cite{AP}, \cite{Sa}, \cite{Ti}) and Sasakian-Einstein manifolds (\cite{BG}). For a homogeneous space the Einstein equation reduces to a system of algebraic equations, which is still very complicated. Even for homogeneous spaces we are far from knowing a full description. For instance, homogeneous Einstein metrics on spheres and projective spaces were classified by Ziller (\cite{Zi2}) and Einstein normal homogeneous manifolds were classified by Wang and Ziller (\cite{WZ}). Every isotropy irreducible space, in particular any irreducible symmetric space (\cite{He},\cite{Ke1}) is an Einstein manifold (\cite{Be}, \cite{Wo}). Recently homogeneous Einstein metrics on homogeneous spaces with exactly two isotropy summands were classified by Dickinson and Kerr (\cite{DK}). Nowadays, it is known that every compact simply-connected homogeneous manifold with dimension less or equal to $11$ admits a homogeneous Einstein metric: a $2$ or $3$-dimensional manifold has constant sectional curvature (\cite{Be}); in dimension $4$, the result was shown by Jensen (\cite{Je1}), and by Alekseevsky, Dotti and Ferraris in dimension 5 (\cite{ADF}); in dimension $6$, the result is due to Nikonorov and Rodionov (\cite{NR2}), and in dimension $7$ it is due to Castellani, Romans and Warner (\cite{CRW}). All the $7$-dimensional homogeneous Einstein manifolds (\cite{Ni}) were obtained by Nikonorov. These results were extended to dimension up to $11$ by B\"{o}hm and Kerr (\cite{BK}). All these results deal with the case of positive Einstein constant. Any homogeneous Riemannian manifold with zero constant is locally flat. There are also many results about invariant Einstein metrics with negative Einstein constant on solvmanifolds.

Einstein homogeneous fibrations have also been the object of study. We recall, for instance, the work of Jensen on principal fibers bundles (\cite{Je2}) and the work of Wang and Ziller on principal torus bundles (\cite{WZ3}).

\li

This paper is devoted to investigation of $G$-invariant Einstein metrics on the total space of $G$-equivariant fibrations. Let $G$ be a compact connected semisimple Lie group and $L\varsubsetneq K\varsubsetneq G$ connected closed non-trivial subgroups of $G$. We consider the fibration

\beq\label{fibDefIntro}M = G/L\rightarrow G/K = N,\textrm{ with fiber } F = K/L\eeq

and investigate the existence of $G$-invariant Einstein metrics on $M$ such that the natural projection $M \ni aL\mapsto aK\in N,\,a\in G,$ is a Riemannian submersion with totally geodesic fibers.

\li

By $\mfg$, $\mfk$ and $\mfl$ we denote the Lie algebras of $G$, $K$ and $L$,
respectively. By $\Phi$ and $\Phi_{\mfk}$ we denote the Killing forms of $G$ and $K$, respectively. We
set $B=-\Phi$ and since $G$ is compact and semisimple, $B$ is positive definite. We denote by $B_{\mathfrak{q}}$ the restriction of $B$ to some subspace $\mathfrak{q}\subset \mfg $.

We consider a $B$-orthogonal
decomposition of $\mfg$ given by

\beq\label{decini}\mfg=\mfl\oplus \mfm,\,\,\, \mfm=\mfp\oplus\mfn,\eeq

where $\mfg=\mfl\oplus\mfm$,
$\mfg=\mfk\oplus\mfn$ and $\mfk=\mfl\oplus\mfp$ are reductive
decompositions for $M$, $N$ and $F$, respectively. An $Ad\,K$-invariant Euclidean product on $\mfn$ induces a
$G$-invariant metric $g_N$ on $N$ and an
$Ad\,L$-invariant Euclidean product on $\mfp$ induces a $G$-invariant metric $g_F$ on $F$. The orthogonal direct sum of the metric $g_F$ on the fiber and the metric $g_N$ on the base space $N$ defines a $G$-invariant
metric $g_M$ on $M$ which projects onto a $G$-invariant metric $g_N$ on
$N$. We recall the following result due to L. B\'{e}rard-Bergery:

\bthm (\cite{BB}, \cite[9 \S H]{Be}) Let $M=G/L\rightarrow G/K=N$ be a $G$-equivariant fibration, for a Lie group $G$, where $K$, $L$ are compact subgroups. Let $g_M$ be the $G$-invariant metric on $M$ given by the orthogonal sum of a $G$-invariant metric $g_N$ on $N$ and a $G$-invariant metric $g_F$ on $F=K/L$. The natural projection $M\ni aL\mapsto aK\in N$ is a Riemannian submersion from $(M,g_M)$ to $(N,g_N)$ with totally geodesic
 fibers.\ethm

Moreover, if $\mfp$ and $\mfn$ do not contain any equivalent
$Ad\,L$-submodules, then any $G$-invariant metric with totally geodesic fibers is constructed in this fashion.

Throughout this paper, a $G$-invariant metric on $M$ such that the natural projection $M\ni aL\mapsto aK\in N$ is a Riemannian submersion with totally geodesic fibers shall be called an \textbf{adapted} metric.

An adapted metric on $M$ shall be denoted by $g_M$; $g_N$ shall denote the projection of
$g_M$ onto the base space $N$ and $g_F$ its restriction to the
fiber $F$. By using the $B$-orthogonal decomposition of $\mfg$ fixed in (\ref{decini}), we denote by $g_{\mfm}$, $g_{\mfn}$ and $g_{\mfp}$ the invariant Euclidean products on $\mfm$, $\mfn$ and $\mfp$, respectively, which determine $g_M$, $g_N$ and $g_F$. Whereas $g_{\mfm}$ and $g_{\mfp}$ are $Ad\,L$-invariant, $g_{\mfn}$ is $Ad\,K$-invariant.

\li

We consider a $B$-orthogonal decomposition $\mfp=\mfp_1\oplus\ldots\oplus\mfp_s$ of $\mfp$ into irreducible $Ad\,L$-submodules and a $B$-orthogonal decomposition $\mfn=\mfn_1\oplus\ldots\oplus\mfn_n$ of $\mfn$ into irreducible $Ad\,K$-submodules.
Throughout we assume the following hypothesis:

\beq\label{submodulesHyp}\bar{l}\mfp_1,\ldots,\mfp_s \textrm{ are pairwise
inequivalent irreducible }Ad\,L\textrm{-submodules;}\\
\mfn_1,\ldots,\mfn_n \textrm{ are
pairwise inequivalent irreducible }Ad\,K\textrm{-submodules;}\\
\mfp \textrm{ and } \mfn \textrm{ do not contain equivalent }
Ad\,L\textrm{-submodules.} \ear\eeq

We remark that the submodules $\mfn_j$ are not required to be $Ad\,L$-irreducible. Under the hypothesis (\ref{submodulesHyp}), according to
Schur's Lemma, any $Ad\,L$-invariant scalar product on
$\mfm=\mfp\oplus\mfn$ which restricts to an $Ad\,K$-invariant
scalar product on $\mfn$ is of the form

\beq\label{mdef}g_{\mfm}=\underbrace{\left(\oplus_{a=1}^s\la_aB_{\mfp_a}\right)}_{g_{\mfp}}\oplus\underbrace{\left(\oplus_{j=1}^n\mu_jB_{\mfn_j}\right)}_{g_{\mfn}},\,\la_a,\,\mu_j>0\eeq

Since an adapted metric $g_M$ on $M$ projects onto a $G$-invariant metric on $N$, $g_M$ is necessarily induced by an Euclidean
product $g_{\mfm}$ of the form (\ref{mdef}).

\li

Throughout this paper any homogeneous fibration $M\rightarrow N$ and any adapted metric $g_M$ on $M$ are as defined above. In Section \ref{sectionRF} we derive formulae for the Ricci curvature for an adapted metric $g_M$
and find some necessary conditions that $g_M$ is Einstein. A main result consists of necessary conditions for existence of an Einstein adapted metric described just in terms of algebraic conditions on the Casimir operators of the isotropy submodules $\mfp_a$ and $\mfn_k$. We recall that if $U$ is a vector subspace of $\mfg$, the Casimir operator of $U$ (with respect to the
Killing form $\Phi$) is the operator

\beq \label{casdef}C_U=\sum_iad_{u_i}ad_{u_i'}\in \mfg\mfl(\mfg),\eeq

where $\{u_i\}_i$ and $\{u_i'\}_i$ are bases of $U$ which are dual
with respect to $\Phi$, i.e., $\Phi(u_i,u_j')=\de_{ij}$. More precisely, we prove the following result:

\bthm \label{cond3}Let $M=G/L\rightarrow G/K$ be a homogeneous fibration. If there exists on $M$ an Einstein adapted
metric, then

(i) there are positive constants $\la_1,\ldots,\la_s$
such that $\sum_{a=1}^s\la_aC_{\mfp_a}$ is scalar on each $\mfn_j$, where $C_{\mfp_a}$ is the Casimir operator of $\mfp_a$;

(ii) if $g_N$ is not a multiple of $B$, then
there are positive constants $\nu_1,\ldots,\nu_n$, which are not all equal,
such that $\sum_{j=1}^n\nu_jC_{\mfn_j}(\mfp)\subset \mfk$, where $C_{\mfn_j}$ is the Casimir operator of $\mfn_j$.
\ethm

In the following sections we focus in some special cases. In Section \ref{sectionBinormal}, we consider the special case when the adapted metric is binormal. A \textbf{binormal} metric on $M$ is a $G$-invariant metric $g_M$ such that its restrictions to the fiber and projection onto the base, $g_F$ and $g_N$, are multiples of the restrictions of the Killing form of $G$, i.e., $g_M$ is defined by the Euclidean product

\beq\label{mdefbi}g_{\mfm}=\la B_{\mfp}\oplus \mu B_{\mfn}.\eeq

A binormal metric is clearly an adapted metric. Binormal metrics are a very natural class of metrics on homogeneous fibrations which generalize normal metrics. We recall that a normal metric on $M$ is defined by the restriction of an $Ad\,G$-invariant positive definite symmetric bilinear map to $\mfm$. Einstein normal homogeneous manifolds were classified by Wang and Ziller (\cite{WZ}) and Rodionov (\cite{Ro1}). A necessary condition for existence of an Einstein normal metric on $M$ is that the Casimir operator of $\mfl$ is scalar on $\mfn$ (\cite{WZ}). Similarly, Theorem \ref{cond3} states that a necessary condition for existence of an Einstein binormal metric is that the Casimir operator of the tangent space to the fibers, $C_{\mfp}$, is scalar on each of the irreducible horizontal submodules $\mfn_j$. If this condition is satisfied, then the following theorem reduces the problem of existence of Einstein binormal metrics to a system of quadratic equations with one variable. Thus, if $C_{\mfp}$ is scalar on each $\mfn_j$, the only difficulty deciding about the existence of an Einstein binormal metric is the computation of the necessary coefficients.

\li

\bthm \label{binormal1}Let $M=G/L\rightarrow G/K$ be a homogeneous fibration.

(i) If the Casimir operator of $\mfp$, $C_{\mfp}$, is not scalar on some
$\mfn_j$, then there are no Einstein binormal metrics on $M$;

(ii) Suppose that $C_{\mfp}$ is scalar on each $\mfn_j$, i.e., $C_{\mfp}\mid_{\mfn_j}=b^jId$. Then
there is a one-to-one correspondence, up to homothety, between Einstein binormal metrics on $M$ and positive solutions of the
following system of quadratic equations on the unknown $X\in\reals$:

\beqar \label{einI1}\de_{ij}^{\mfk}(1-X)=\de_{ij}^{\mfl}, \textrm{ if
}n>1, \\
\label{einI2}(2\de_{ab}^{\mfl}+\de_{ab}^{\mfk})X^2=\de_{ab}^{\mfk},
\textrm{ if } s>1, \\
\label{einI3}\left(\ga_a+2c_{\mfl,a}\right)X^2-\left(1+2c_{\mfk,j}\right)X+(1-\ga_a+2b^j)=0,\eeqar

for $a,b=1,\ldots,s$ and $i,j=1,\ldots,n$. Here $c_{\mfl,a}$
is the eigenvalue of $C_{\mfl}$ on $\mfp_a$, $\ga_a$ is the constant determined
by

$$\Phi_{\mfk}\mid_{\mfp_a\times\mfp_a}=\ga_a\Phi\mid_{\mfp_a\times\mfp_a},$$

$c_{\mfk,j}$ is the eigenvalue of $C_{\mfk}$ on $\mfn_j$ and the
$\de$'s are the differences $\de_{ij}^{\mfk}=c_{\mfk,i}-c_{\mfk,j}$, $\de_{ij}^{\mfl}=c_{\mfl,i}-c_{\mfl,j}$, $\de_{ab}^{\mfk}=\ga_a-\ga_b$ and $\de_{ab}^{\mfl}=c_{\mfl,a}-c_{\mfl,b}$.

If such a positive solution $X$ exists, then Einstein binormal metrics are, up to
homothety, given by
$$g_{\mfm}=B_{\mfp}\oplus XB_{\mfn}.$$
\ethm

In particular, this allows us to characterize Einstein adapted metrics on fibrations such that the fiber and base spaces are isotropy irreducible spaces. We remark that the fact that the base space $N$ is isotropy irreducible does not imply that $M$ has only two irreducible isotropy subspaces, since the horizontal subspace $\mfn$ is $Ad\,K$-irreducible but is not in general $Ad\,L$-irreducible. The particular case of existence of $G$-invariant Einstein metrics when $M$ has exactly two irreducible isotropy subspaces was studied by Wang and Ziller, under some assumptions, in \cite{WZ2} and more recently, in full generality, by Dickinson and Kerr in \cite{DK}, who classified all such metrics.

\bcor\label{bincor1} Let $M=G/L\rightarrow G/K=N$ be a homogeneous fibration. Suppose that the fiber $F=K/L$ and the base space $N$ are isotropy
irreducible spaces and $dim\,F>1$. There exists on $M$ an
Einstein adapted metric if and only if the following two conditions are satisfied:

(i) $C_{\mfp}$ is scalar on
$\mfn$;

(ii) $\triangle \geq 0$, where
$$\triangle=(1+2c_{\mfk,\mfn})^2-4(\ga+2c_{\mfl,\mfp})(1-\ga+2b),$$

$c_{\mfk,\mfn}$ is the eigenvalue of $C_{\mfk}$ on $\mfn$,
$c_{\mfl,\mfp}$ is the eigenvalue of $C_{\mfl}$ on $\mfp$, $b$ is
the eigenvalue of $C_{\mfp}$ on $\mfn$ and $\ga$ is such that
$\Phi_{\mfk}\mid_{\mfp\times\mfp}=\ga \Phi\mid_{\mfp\times\mfp}$.

If these conditions are satisfied, then Einstein adapted
metrics are, up to homothety, given by

$$g_{\mfm}=B_{\mfp}\oplus XB_{\mfn}, \textrm{
where
}X=\frac{1+2c_{\mfk,\mfn}\pm\sqrt{\triangle}}{2(\ga+2c_{\mfl,\mfp})}.$$\ecor

In the case when $F$ is 1-dimensional, the fibration $M\rightarrow
N$ is a principal circle bundle, since $F$ is an abelian compact connected
$1$-dimensional group. We recall that Einstein metrics on principal fiber bundles have been widely studied (\cite{Je2},\cite{WZ3}) and, in particular, homogeneous Einstein metrics on circle bundles were classified McKenzie Y. Wang and Wolfgang Ziller in \cite{WZ3}. From Theorem \ref{binormal1}, We obtain the following well known result obtained previously in \cite{WZ3}.

\bcor\label{bincor2} Let $M=G/L\rightarrow G/K=N$ be a homogeneous fibration. Suppose that $N$ is isotropy irreducible and the fiber $F=K/L$ is isomorphic to the circle group. There exists on $M$ exactly
one $G$-invariant Einstein metric, up to homothety, given by

$$g_{\mfm}=B_{\mfp}\oplus XB_{\mfn},\,where\,X=\frac{2+m}{m(1+2c_{\mfk,\mfn})},$$

$c_{\mfk,\mfn}$ is the eigenvalue of $C_{\mfk}$ on $\mfn$
and $m=dim\,N$. \ecor

Also interesting necessary conditions are found if the fiber is not isotropy irreducible.

\bcor \label{binormal3} Let $M=G/L\rightarrow G/K$ be a homogeneous fibration. Suppose the fiber $F=K/L$ is not isotropy irreducible and there exists a constant $\al$ such that

\beq \Phi\circ C_{\mfl}\mid_{\mfp\times\mfp}=\al\Phi_{\mfk}\mid_{\mfp\times\mfp}.\eeq

For $a=1,\ldots,s$, let $\ga_a$ be the constant defined by
$\Phi_{\mfk}\mid_{\mfp_a\times\mfp_a}=\ga_a\Phi\mid_{\mfp_a\times\mfp_a}$. If $\ga_a\neq \ga_b$, for some $a,\,b$, then there exists an
Einstein binormal metric on $M$ if and only if $C_{\mfp}$ is scalar on each $\mfn_j$ and

\beq\label{caskjlj}c_{\mfl,j}=\left(1-\frac{1}{\sqrt{2\al+1}}\right)\left(c_{\mfk,j}+\frac{1}{2}
\right),\,j=1,\ldots,n,\eeq

where $c_{\mfl,j}$ and
$c_{\mfk,j}$ are the eigenvalues of $C_{\mfl}$ and $C_{\mfk}$,
respectively, on $\mfn_j$. In this case, there is a unique binormal
Einstein metric, up to homothety, given by

$$g_{\mfm}=B_{\mfp}\oplus
\frac{1}{\sqrt{2\al+1}}B_{\mfn}$$

and, furthermore, the number $\sqrt{2\al+1}$ is rational. \ecor

\li

Interesting applications arise from the result above, for instance, when the fiber $F=K/L$ is a symmetric space, since, in this case, $\al=1/2$ and clearly $\sqrt{2\al+1}$ is not rational. This implies that if an Einstein binormal metric exists, then the Casimir operator of $\mfk$ must be scalar on $\mfp$. Though, Einstein adapted metrics on homogeneous fibrations with symmetric fiber are out of scope of this paper, the reader is invited to find such applications in \cite{Fa}.

A natural question is to determine an Einstein adapted  metric whose restriction to the fiber and projection to the base space are also Einstein metrics. This problem has been approached by several authors, see for instance, results by Berard-Bergery, Matsuzawa and Wang and Ziller in \cite{Be} and \cite{WZ}. In Section \ref{einsteinFB} we prove the following two results:

\bthm \label{bfein} Let $g_M$ be an Einstein adapted metric on the homogeneous fibration $M=G/L\rightarrow G/K=N$ defined by the Euclidean product

$$g_{\mfm}=\left(\oplus_{a=1}^s\la_aB_{\mfp_a}\right)\oplus\left(\oplus_{k=1}^n\mu_kB_{\mfn_k}\right).$$

If $g_N$ and $g_F$ are also Einstein, then

\beqar\frac{\mu_j}{\mu_k}=\left(\frac{b^j}{b^k}\right)^{\frac{1}{2}},\,j,k=1,\ldots,n,\\
\frac{\la_a}{\la_b}= \sum_{j=1}^n\frac{C_{\mfn_j,b}}{b_j}\Big/\sum_{j=1}^n\frac{C_{\mfn_j,a}}{b_j},\,a,b=1,\ldots,s,\eeqar

where $b^j$ is the eigenvalue of  $\sum_{a=1}^s\la_aC_{\mfp_a}$ on $\mfn_j$ and $c_{\mfn_j,a}$ is defined by $$\Phi(C_{\mfn_j}\cdot,\cdot)\mid_{\mfp_a\times\mfp_a}=c_{\mfn_j,a}\Phi\mid_{\mfp_a\times\mfp_a}.$$

In particular, there exists at most one $G$-invariant metric $g_N$ on $N$ and one $K$-invariant metric $g_F$ on $F$ such that $g_M$ is Einstein.

\ethm

\li

\bthm \label{einrel} Let $g_M$ be an adapted metric on the homogeneous fibration $M=G/L\rightarrow G/K=N$ defined by the Euclidean product

$$g_{\mfm}=\left(\oplus_{a=1}^s\la_aB_{\mfp_a}\right)\oplus\left(\oplus_{k=1}^n\mu_kB_{\mfn_k}\right).$$

Suppose that $g_M$, $g_N$ and $g_F$ are
Einstein and let $E$, $E_F$ and $E_N$ be the corresponding
Einstein constants. If $E\neq E_N$, then

$$\mu_j=\left(\frac{b^j}{2(E_N-E)}\right)^{\frac{1}{2}},$$

$$\la_a=2\frac{E-E_F}{E_N-E}\left(\sum_{j=1}^n\frac{C_{\mfn_j,a}}{b^j}\right)^{-1}.$$

where $b^j$ is the eigenvalue of the operator
$\sum_{a=1}^s\la_aC_{\mfp_a}$ on $\mfn_j$ and $c_{\mfn_j,a}$ is defined by $\Phi(C_{\mfn_j}\cdot,\cdot)\mid_{\mfp_a\times\mfp_a}=c_{\mfn_j,a}\Phi\mid_{\mfp_a\times\mfp_a}$.

\ethm

Finally, in Section \ref{kow}, we prove the existence of a non-standard Einstein metric on the $n$-symmetric spaces $M=\frac{G_0^n}{\triangle^nG_0}$, for a compact connected simple Lie group $G_0$, where $\triangle^nG_0$ is the diagonal subgroup. We recall that Wang and Ziller have shown that the standard metric on these spaces is Einstein (\cite{WZ2}, \cite{Ro1}). To investigate the existence of a non-standard Einstein $G$-invariant metric on $M$, we consider the fibration

$$\frac{G_0^n}{\triangle^nG_0}\rightarrow
\frac{G_0^p}{\triangle^pG_0}\times \frac{G_0^q}{\triangle^qG_0},\textrm{ with fiber } F=\frac{\triangle^p G_0\times \triangle^q G_0}{\triangle^n G_0},$$

for $n=p+q$. The vertical isotropy subspace $\mfp$ is $Ad\,L$-irreducible. However, the horizontal isotropy subspace $\mfn$ is not $Ad\,K$-irreducible. We consider a decomposition $\mfn=\mfn_1\oplus\mfn_2$,

\beqar\label{decnsKowalski}\mfn_1=\{(X_1,\ldots,X_p,0,\ldots,0):X_j\in\mfg_0,\sum X_j=0
\}\subset \mfg_0^p\times 0_q\\ \nonumber
\mfn_2=\{(0,\ldots,0,X_1,\ldots,X_q):X_j\in\mfg_0,\sum X_j=0
\}\subset 0_p\times\mfg_0^q,\eeqar

where $\mfg_0$ is the Lie algebra of $G_0$. We remark that this decomposition of $\mfn$ is not unique and $\mfn_1$, $\mfn_2$ are not $Ad\,K$-irreducible either. Hence, the hypothesis (\ref{submodulesHyp}) is not satisfied. We can still consider an adapted metric $g_M$ defined by a scalar product of the form

\beq\label{defmKow}g_{\mfm}=\la B_{\mfp}\oplus \mu_1 B_{\mfn_1}\oplus \mu_2 B_{\mfn_2},\eeq

but we should keep in mind that other adapted metrics might exists which are not of this form. We prove the following result:

\bthm\label{genknss} Let $G_0$ be a compact connected simple Lie group and consider the fibration
$$M=\frac{G_0^n}{\triangle^nG_0}\rightarrow
\frac{G_0^p}{\triangle^pG_0}\times \frac{G_0^q}{\triangle^qG_0},$$

where $p+q=n$ and $2\leq p\leq q\leq n-2$. For $n=4$, the only Einstein adapted metric is the standard metric. If $n>4$, there exist on $M$ at least one non-standard Einstein adapted metric of the form (\ref{defmKow}). This non-standard Einstein adapted metric is binormal if and only if $p=q$. Furthermore, its projection onto the base space is also Einstein if and only if $p=q$.

\ethm

\section{The Ricci Curvature}\label{ricciCurv}

Let $M=G/L$ be a homogeneous manifold of a connected Lie group $G$, where $L$ is a compact subgroup. Let $\mfg=\mfl\oplus\mfm$ be a reductive decomposition of $M$. In this section we describe the Ricci curvature, $Ric$, of the $G$-invariant metric $g_M$ on $M$ associated to the
$Ad\,L$-invariant scalar product $<,>$ on $\mfm$. For $X\in\mfg$, let $P_X$ and $T_X$ be the endomorphisms of $\mfm$ defined by

\beqar P_XY=[X,Y]_{\mfm}, \, Y\in\mfm\\
\label{opT}<T_XY,Z>=<X,P_YZ>,\,Y,Z\in\mfm \eeqar

where the subscript $\mfm$ denotes projection onto $\mfm$. For $X\in\mfg$, The \textit{Nomizu operator} $L_X$ of the scalar product $<,>$ (\cite{No}, \cite{NK}) is defined by

\beq L_XY=-\nabla_{Y}X^*,\,Y\in\mfm\eeq

where $\nabla$ is the Riemannian connection of $g_M$ and $X^*$ is the Killing vector field generated by $X$. We have

\beq \label{no1}L_XY=\frac{1}{2}P_XY+U(X,Y),\eeq

where $U:\mfm\times\mfm\rightarrow \mfm$ is the operator

\beq\label{opU}U(X,Y)=-\frac{1}{2}(T_XY+T_YX),\,X,Y\in\mfm.\eeq

Moreover, $L_X$ is skew-symmetric with respect to $<,>$ and $L_XY-L_YX=P_XY$, $Y\in\mfm$. The metric $g_M$ is called \textit{naturally reductive} if $U=0$. The curvature tensor of $g_M$ at the point $o=eL$ can be written as

\beq
\label{curv}R(X,Y) = [L_X,L_Y]-L_{[X,Y]},\,\,\,X,\,Y\in \mfm=T_oM.\eeq

The sectional curvature $K$ of $g_M$ is defined by

\beq K(Z,X)=<R(Z,X)X,Z>,\eeq

for every $X,Z\in\mfm$ orthonormal with respect to $<,>$. The Ricci curvature of $g_M$ is determined by

\beq Ric(X,X)=\sum_iK(Z_i,X),\,X\in\mfm\eeq

where $(Z_i)_i$ is an orthonormal basis of $\mfm$ with respect to $<,>$. The metric $g_M$ is said to be an Einstein metric if there exists a constant $E$ such that $Ric=Eg_M$. Below we describe the Ricci curvature of the $G$-invariant metric $g_M$ by using the endomorphisms described above. The proof is out of the scope of this paper and can be found in detail in \cite[\S 1.1]{Fa}.

\blem\label{ric2} \cite[\S 1.1]{Fa} Let $X,Y\in\mfm$.

$$Ric(X,Y)=-\frac{1}{4}tr(2P_X^*P_Y+T_XT_Y)-\frac{1}{2}\Phi(X,Y)+tr(P_{U(X,Y)}).$$

\elem

\brem \label{ric2rem}If there exists on $\mfm$ a non-degenerate $Ad\,L$-invariant symmetric bilinear form $\be$, then $tr(P_{U(X,Y)})=0$, for all $X,Y\in\mfm$. Indeed, if such a
bilinear form exists, $tr\,P_a=0$, for every $a\in\mfm$. Let
$\{w_i\}_i$ and $\{w_i'\}_i$ be bases of $\mfm$ dual with respect
to $\be$, i.e., $\be(w_i,w_j')=\de_{ij}$. Then, for every
$a\in\mfm$,

$$\be(P_aw_i,w_i')= \be([a,w_i]_{\mfm},w_i')=-\be(w_i,[a,w_i']_{\mfm})= -\be(P_aw_i',w_i).$$

Hence, $tr(P_a)=0$. Also, if the metric $g_M$ on $M$ is naturally reductive, then $P_{U(X,Y)}=0$, for all $X,Y\in\mfm$, since, in this case, $U$ is identically zero.

$\diamond$

\erem

\bdfn \label{Casimirdef}Let $\be$ be a non-degenerate $Ad\,G$-invariant symmetric bilinear form on $\mfg$. Let $U$ be an $Ad\,L$-invariant vector
subspace of $\mfg$ such that the restriction of $\be$ to $U$ is
non-degenerate. The Casimir operator of $U$ with respect to $\be$ is the operator

$$C_U=\sum_iad_{u_i}ad_{u_i'},$$

where $\{u_i\}_i$ and $\{u_i'\}_i$ are bases of $U$ which are dual
with respect to $\be$, i.e., $\be(u_i,u_j')=\de_{ij}$.\edfn

A Casimir operator is independent of the choice of dual basis. Moreover, it is an $Ad\,L$-invariant linear map
and thus it is scalar on any irreducible $Ad\,L$-module. In
particular, if $\mfg$ is simple, the only non-degenerate $Ad\,G$-invariant symmetric bilinear map on $\mfg$, up to scalar factor, is the Killing form $\Phi$ and $C_{\mfg}=Id$.

\bdfn \label{tracesdef}Let $U$, $V$ be $Ad\,L$-invariant vector
subspaces of $\mfg$. We define a bilinear map $Q_{UV}:\mfm\times\mfm
\rightarrow \reals$ by

$$Q_{UV}(X,Y)=tr([X,[Y,\cdot]_V]_U), \,X,Y\in\mfm,$$

where the subscripts $U$ and $V$ denote the projections onto $U$
and $V$, respectively. \edfn

\li

\blem \label{traces1} Let $\be$ be a non-degenerate $Ad\,G$-invariant symmetric bilinear form on $\mfg$. Let $U$,
$V$ be $Ad\,L$-invariant vector subspaces of $\mfg$ such that the
restrictions of $\be$ to $U$ and $V$ are both non-degenerate. Then

(i) $Q_{UV}=Q_{VU}$ and $Q_{UV}$ is an $Ad\,L$-invariant symmetric bilinear map.
Hence, if $W\subset\mfg$ is any irreducible $Ad\,L$-submodule,
then $Q_{UV}\mid_{W\times W}$ is a multiple of
 $\be\mid_{W\times W}$.

(ii) if $ad_XU\subset V$ or $ad_YU\subset V$, then
$Q_{UV}(X,Y)=\be(C_UX,Y)=\be(X,C_UY)$, for $X,Y\in\mfm$, where $C_U$ and $C_V$ are the Casimir operators of $U$ and $V$, respectively, with respect to $\be$;

(iii) if $ad_XV\perp U$ or $ad_YV\perp U$, then $Q_{UV}(X,Y)=0$, for $X,Y\in\mfm$;

(iv) if $ad_Xad_YU\perp U$ or $ad_Yad_XU\perp U$, then
$Q_{UV}(X,Y)=0$, for $X,Y\in\mfm$.

\elem

\bproof Since $\be$ is non-degenerate and $\be\mid_{U\times U}$, $\be\mid_{V\times V}$ are non-degenerate, we consider the
orthogonal complements $U^{\perp}$ and $V^{\perp}$ of $U$ and $V$,
respectively, in $\mfg$ with respect to $\be$. Also, we consider bases $\{w_i\}_i$ and $\{w_i'\}_i$ of $U$
which are dual with respect to $\be$. Moreover, since $\be$ is $Ad\,L$-invariant, it is associative. Let $X,Y\in\mfm$ and $g\in L$.

$$\bar{rl}\be([X,[Y,w_i]_V]_U,w_i')= & \be([X,[Y,w_i]_V],w_i')\\ \\

= & -\be([Y,w_i]_V,[X,w_i'])\\ \\

= & -\be([Y,w_i],[X,w_i']_V)\\ \\

= & \be(w_i,[Y,[X,w_i']_V])\\ \\

= & \be(w_i,[Y,[X,w_i']_V]_U). \ear$$

Therefore, $tr([X,[Y,\cdot]_V]_U)=tr([Y,[X,\cdot]_V]_U)$ and thus
$Q_{UV}(X,Y)=Q_{UV}(Y,X)$. So $Q_{UV}$ is symmetric. To show the $Ad\,L$-invariance of $Q_{UV}$ we note that since $V$
and $V^{\perp}$ are $Ad\,L$-invariant subspaces and $\mfg=V\oplus
V^{\perp}$, the projections on $V$ and $V^{\perp}$ are also
$Ad\,L$-invariant linear maps.

$$\bar{rl}\be([Ad_gX,[Ad_gY,w_i]_V]_U,w_i')= & \be([Ad_gX,[Ad_gY,w_i]_V],w_i')\\ \\

= & \be(Ad_{g^{-1}}[Ad_gX,[Ad_gY,w_i]_V],Ad_{g^{-1}}w_i')\\ \\

= & \be([X,Ad_{g^{-1}}[Ad_gY,w_i]_V],Ad_{g^{-1}}w_i')\\ \\

= & \be([X,[Y,Ad_{g^{-1}}w_i]_V],Ad_{g^{-1}}w_i')\\ \\

= & \be([X,[Y,Ad_{g^{-1}}w_i]_V]_U,Ad_{g^{-1}}w_i').\ear$$

Since $\be$ is $Ad\,G$-invariant, if $\{w_i\}_i$ and $\{w'_i\}_i$ are dual bases of $U$ with respect
to $\be$, then $\{Ad_{g^{-1}}w_i\}_i$ and $\{Ad_{g^{-1}}w'_i\}_i$
are still dual bases as well. So by the above we conclude that $tr([Ad_gX,[Ad_gY,\cdot]_V]_U)=tr([X,[Y,\cdot]_V]_U)$ and thus
$Q_{UV}$ is $Ad\,L$-invariant.

Let $Z\in\mfm$ and set $A_Z=\left(ad_Z\mid_U\right)_V$ and
$B_Z=\left(ad_Z\mid_V\right)_U$. We have

$$Q_{VU}(X,Y)=tr(A_XB_Y)=tr(B_YA_X)=Q_{UV}(Y,X).$$

Hence, by symmetry of $Q_{UV}$, we conclude that
$Q_{VU}(X,Y)=Q_{UV}(Y,X)=Q_{UV}(X,Y)$, for every $X,Y\in\mfm$.
Therefore, $Q_{UV}=Q_{VU}$.

To show (ii) let $C_U=\sum_iad_{w_i}ad_{w'_i}$ be the Casimir
operator of $U$ with respect to $\be$. Since $Q_{UV}=Q_{VU}$ it suffices to suppose that
$ad_YU\subset V$. If $ad_YU\subset V$, then
$$Q_{UV}(X,Y)=tr([X,[Y,\cdot]]_U)=tr(ad_Xad_Y\mid_U).$$

Since $\be([X,[Y,w_i]]_U,w'_i) = \be([X,[Y,w_i]],w'_i)= \be(Y,[w_i,[w'_i,X]])$,  we have

$$Q_{UV}(X,Y)=\sum_i\be(Y,[w_i,[w'_i,X]])=\be(Y,C_UX).$$

By symmetry of $Q_{UV}$ we also get $Q_{UV}(X,Y)=\be(X,C_UY)$.

If $ad_XV\perp U$, then, for every $w,w'\in U$,
$\be([X,[Y,w]_V],w')=0$ and thus $Q_{UV}(X,Y)=0$, for every $Y\in
\mfm$. By symmetry, the same conclusion holds if $ad_YV\perp U$. This shows (iii).

Finally, to prove (iv), if $ad_Xad_YU\perp U$, then, for every $w,w'\in U$,
$\be([X,[Y,w]],w')=0$ and thus $\be([X,[Y,w]_V]_U,w')=0$. Hence
$Q_{UV}(X,Y)=0$. If $ad_Yad_XU\perp U$, then $Q_{UV}(X,Y)=0$ by
symmetry.

$\Box$ \eproof

The following theorem describes the Ricci curvature of an invariant metric on $M=G/L$, when the Lie algebra admits a non-degenerate $Ad\,G$-invariant symmetric bilinear form.

\bthm \label{ric1} Let $M=G/L$ be a homogeneous manifold and let $\be$ be a non-degenerate $Ad\,G$-invariant symmetric bilinear form on the Lie algebra $\mfg$. Let $g_M$ be the $G$-invariant
pseudo-Riemannian metric on $M$ induced by the scalar product of
the form

\beq\label{mdef1}<,>=\oplus_{j=1}^m\nu_j \be\mid_{\mfm_j\times\mfm_j},\,\nu_j>0,\eeq

where $\mfm=\mfm_1\oplus\ldots\oplus\mfm_m$ is a $\be$-orthogonal decomposition of $\mfm$. For $X\in\mfm_a$, $Y\in\mfm_b$, the Ricci curvature of $g_M$ is given by

$$Ric(X,Y)=\frac{1}{2}\sum_{j,k=1}^m\left(\frac{\nu_k}{\nu_j}-\frac{\nu_a\nu_b}{2\nu_k\nu_j}\right)Q_{\mfm_j\mfm_k}(X,Y)-\frac{1}{2}\Phi(X,Y).$$

\ethm

\bproof First we note that
$\be\mid_{\mfm_j\times\mfm_j}$ is in fact non-degenerate. Let
$X\in\mfm_a$ and $Y\in\mfm_b$. We apply the formula given in Lemma \ref{ric2}. According to Remark \ref{ric2rem}, we have $tr(P_{U(X,Y)})=0$. Let $j=1,\ldots,m$ and let $\{w_i\}_i$ and $\{w'_i\}_i$ be dual
bases for $\mfm_j$ with respect to $\be$.

$$\bar{rl}<T_XT_Yw_i,w_i'>= & <X,[T_Yw_i,w'_i]_{\mfm}>\\ \\

= & \nu_a \be(X,[T_Yw_i,w'_i])\\ \\

= & -\nu_a \be(T_Yw_i,[X,w'_i])\\ \\

= & -\nu_a \sum_{k=1}^m \be(T_Yw_i,[X,w'_i]_{\mfm_k})\\ \\

= & -\nu_a \sum_{k=1}^m\nu_k^{-1}<T_Yw_i,[X,w'_i]_{\mfm_k}>\\ \\

= & -\nu_a \sum_{k=1}^m\nu_k^{-1}<Y,[w_i,[X,w'_i]_{\mfm_k}]_{\mfm}>\\ \\
= & -\nu_a\nu_b \sum_{k=1}^m\nu_k^{-1}\be([Y,w_i],[X,w'_i]_{\mfm_k})\\ \\

= & -\nu_a\nu_b \sum_{k=1}^m\nu_k^{-1}\be([Y,w_i]_{\mfm_k},[X,w'_i])\\ \\

= & \nu_a\nu_b \sum_{k=1}^m\nu_k^{-1}\be(w'_i,[X,[Y,w_i]_{\mfm_k}])\\ \\

= & \nu_a\nu_b \sum_{k=1}^m\nu_k^{-1}\be(w'_i,[X,[Y,w_i]_{\mfm_k}]_{\mfm_j})\\ \\

= & \nu_a\nu_b
\sum_{k=1}^m(\nu_k\nu_j)^{-1}<w'_i,[X,[Y,w_i]_{\mfm_k}]_{\mfm_j}>.\ear$$

This implies that

$$tr(T_XT_Y\mid_{\mfm_j})=\nu_a\nu_b
\sum_{k=1}^m\dfrac{1}{\nu_k\nu_j}tr([X,[Y,\cdot]_{\mfm_k}]_{\mfm_j})=\nu_a\nu_b
\sum_{k=1}^m\dfrac{1}{\nu_k\nu_j}Q_{\mfm_j\mfm_k}(X,Y)$$

and thus $tr(T_XT_Y)=\nu_a\nu_b
\sum_{j,k=1}^m\frac{1}{\nu_k\nu_j}Q_{\mfm_j\mfm_k}(X,Y)$.

\li

$$\bar{rl}<P_X^*P_Yw_i,w'_i>= & <P_Yw_i,P_Xw'_i>\\ \\

= & \sum_{k=1}^m<[Y,w_i]_{\mfm_k},[X,w'_i]_{\mfm_k}>\\ \\

= & \sum_{k=1}^m\nu_k\be([Y,w_i]_{\mfm_k},[X,w'_i])\\ \\

= & -\sum_{k=1}^m\nu_k\be(w'_i,[X,[Y,w_i]_{\mfm_k}])\\ \\

= & -\sum_{k=1}^m\nu_k\be(w'_i,[X,[Y,w_i]_{\mfm_k}]_{\mfm_j})\\ \\

= &
-\sum_{k=1}^m\nu_k\nu_j^{-1}<w'_i,[X,[Y,w_i]_{\mfm_k}]_{\mfm_j}>.\ear$$

Then
$$tr(P_X^*P_Y\mid_{\mfm_j})=-\sum_{k=1}^m\dfrac{\nu_k}{\nu_j}tr([X,[Y,\cdot]_{\mfm_k}]_{\mfm_j})=-
\sum_{k=1}^m\dfrac{\nu_k}{\nu_j}Q_{\mfm_j\mfm_k}(X,Y)$$

and thus we get $tr(P_X^*P_Y)=-\sum_{j,k=1}^m\frac{\nu_k}{\nu_j}Q_{\mfm_j\mfm_k}(X,Y)$.

By using Lemma \ref{ric2} we finally obtain the required
expression for $Ric(X,Y)$.

$\Box$ \eproof

We recall that a metric $g_M$ is said to be \textit{normal} if it is the restriction of a non-degenerate $Ad\,L$-invariant
symmetric bilinear form on $\mfm$. The formula below for the Ricci curvature of a normal metric was first found by Wang and Ziller in \cite{WZ}, and can be deduced from Theorem \ref{ric1}. From Corollary \ref{ric1cor1}, it is clear that a necessary and sufficient condition for a normal metric to be Einstein is that the Casimir operator of $\mfl$ is scalar on the isotropy space $\mfm$. For instance, this condition holds if $\mfm$ is irreducible.  simply-connected non-strongly isotropy irreducible homogeneous spaces which admit a normal Einstein metric were classified by Wang and Ziller in \cite{WZ}, when $G$ is a compact connected simple group. Also, more generally, simply-connected compact standard homogeneous manifolds were studied by E.D. Rodionov in \cite{Ro4}.

\bcor \label{ric1cor1}Let $\be$ be a non-degenerate $Ad\,G$-invariant symmetric bilinear form on $\mfg$ and $g_M$ the normal metric on $M$ defined by the restriction of $\be$ to $\mfm$. The Ricci curvature of $g_M$ is given by $Ric(\mfm_a,\mfm_b)=0$, if $a\neq b$, and for $X\in\mfm_a$,

$$Ric(X,X)=-\frac{1}{4}\Phi(X,X)-\frac{1}{2}\be(C_{\mfl}X,X),$$

where $C_{\mfl}$ is the Casimir operator of $\mfl$ with respect to
$\be$. \ecor

\bproof Let $C_{\mfg}$, $C_{\mfl}$ and $C_{\mfm}$ be the Casimir operators of $\mfg$, $\mfl$ and $\mfm$ with respect to $\be$. We remark that the Killing form of $\mfg$ may not be non-degenerate. Since $g_M$ is defined by the restriction of $\be$, in Theorem \ref{ric1} we can take $\nu_1=\ldots=\nu_m=1$
and obtain the following. Let $X\in\mfm_a$ and $Y\in\mfm_b$.

$$\bar{rl}Ric(X,Y)= &
\frac{1}{4}\sum_{j,k=1}^mQ_{\mfm_j\mfm_k}(X,Y)-\frac{1}{2}\Phi(X,Y)\\ \\

= & \frac{1}{4}Q_{\mfm\mfm}(X,Y)-\frac{1}{2}\Phi(X,Y)\\ \\

= & \frac{1}{4}Q_{\mfm\mfg}(X,Y)-\frac{1}{4}Q_{\mfm\mfl}(X,Y)-\frac{1}{2}\Phi(X,Y)\\ \\

= & \frac{1}{4}\be(C_{\mfm}X,Y)-\frac{1}{4}\be(C_{\mfl}X,Y)-\frac{1}{2}\Phi(X,Y)\\ \\

= & -\frac{1}{4}\Phi(X,Y)-\frac{1}{2}\be(C_{\mfl}X,Y).\ear$$

Since $C_{\mfl}(\mfm_a)\subset\mfm_a$, it is clear that
$Ric(X,Y)=0$ if $a\neq b$ and $Ric$ is well determined by elements
$Ric(X,X)$ with $X\in\mfm_a$.

$\Box$\eproof

\li

We obtain a similar formula to that of Corollary \ref{ric1cor1}, in the case when the submodules $\mfm_1,\ldots,\mfm_m$ pairwise commute. The proof is similar.

\bcor \label{ric1cor2} Let $\be$ be a non-degenerate $Ad\,G$-invariant symmetric bilinear form on $\mfg$ and $g_M$ on $M$ defined by (\ref{mdef1}). Suppose that $[\mfm_a,\mfm_b]=0$, if $a\neq b$. Then
$Ric(\mfm_a,\mfm_b)=0$, for $a\neq b$, and for $X\in\mfm_a$,

$$Ric(X,X)=-\frac{1}{4}\Phi(X,X)-\frac{1}{2}\be(C_{\mfl}X,X),$$

where $C_{\mfl}$ is the Casimir operator of $\mfl$ with respect to
the $\be$. \ecor

\li

\section{The Ricci Curvature of an Adapted Metric}\label{sectionRF}

In this section we obtain the Ricci curvature of an adapted
metric $g_M$ on the total space of a
homogeneous fibration as in (\ref{fibDefIntro}). Let  $$M=G/L\rightarrow G/K=N,$$

with fiber $F=K/L,$  for a compact connected semisimple Lie group $G$ and $$\mfg=\mfl\oplus \mfm=\mfl\oplus\mfp\oplus\mfn$$ an associated reductive decomposition. We use the notation convention from Section \ref{intro}. We recall that an adapted metric $g_M$ on $M$ is induced by an $Ad\,L$-invariant Euclidean product $g_{\mfm}$ given by (\ref{mdef}), i.e.,
$$g_{\mfm}=\left(\oplus_{a=1}^s\la_aB_{\mfp_a}\right)\oplus\left(\oplus_{k=1}^n\mu_kB_{\mfn_k}\right).$$

\li

All the Casimir operators $C_{\mfk}$, $C_{\mfp_a}$ and $C_{\mfn_j}$ are with respect to the Killing form $\Phi$ (see Definition \ref{casdef}). Since $\Phi(C_{\mfk}\cdot, \cdot)$ and $\Phi(C_{\mfn_j}\cdot, \cdot)$ are $Ad\,L$-invariant symmetric bilinear maps and $\mfp_a$ is $Ad\,L$-irreducible, there are constants $\ga_a$ and $c_{\mfn_j,a}$ such that

\beqar\label{gas}\Phi(C_{\mfk}\cdot,\cdot)\mid_{\mfp_a\times\mfp_a}=\ga_a\Phi\mid_{\mfp_a\times\mfp_a}\\
\label{cnjs}\Phi(C_{\mfn_j}\cdot,\cdot)\mid_{\mfp_a\times\mfp_a}=c_{\mfn_j,a}\Phi\mid_{\mfp_a\times\mfp_a}.\eeqar

\li

In the following Lemma we use the bilinear form $Q_{UV}$ from Definition \ref{tracesdef}.

\blem \label{tracespn}Let $X\in\mfp$ and $Y\in\mfm$.

(i) $Q_{\mfn_j\mfp_a}(X,Y)=Q_{\mfp_a\mfn_j}(X,Y)=0$;

(ii) $Q_{\mfn_i\mfn_j}(X,Y)=0$, if $i\neq j$;

(iii) $Q_{\mfn_j\mfn_j}(X,Y)=\Phi(C_{\mfn_j}X,Y)$.

\li

Let $X'\in\mfn_k$ and $Y'\in\mfm$.

(iv) $Q_{\mfn_j\mfp_a}(X',Y')=Q_{\mfp_a\mfn_j}(X',Y')=0$, if $j\neq
k$;

(v)
$Q_{\mfp_a\mfn_k}(X',Y')=Q_{\mfn_k\mfp_a}(X',Y')=\Phi(C_{\mfp_a}X',Y')$;

(vi) $Q_{\mfp_b\mfp_a}(X',Y')=0$.
\elem

\bproof Let $X\in\mfp$ and $Y\in\mfm$. Since $ad_X\mfp\subset\mfk\perp\mfn$ we have
$Q_{\mfn_j\mfp_a}(X,Y)=0$, from Lemma \ref{traces1}. From Lemma
\ref{traces1},
$Q_{\mfp_a\mfn_j}(X,Y)=Q_{\mfn_j\mfp_a}(X,Y)=0$. As $ad_X\mfn_j\subset\mfn_j$, we have
$Q_{\mfn_j\mfn_j}(X,Y)=\Phi(C_{\mfn_j}X,Y)$. Moreover, since
$\mfn_j\perp\mfn_i$, for every $i\neq j$, we also conclude that
$Q_{\mfn_i\mfn_j}(X,Y)=0$, if $i\neq j$.

\li

Let $X'\in\mfn_k$ and $Y'\in\mfm$. We have $ad_X'\mfp_a\subset\mfn_k\perp\mfp,\mfn_j$, for
$j\neq k$. Thus, $Q_{\mfn_j\mfp_a}(X',Y')=0$, if $j\neq k$ and
$Q_{\mfp_j\mfp_a}(X',Y')=0$, from Lemma
\ref{traces1}. Also from $ad_X'\mfp_a\subset\mfn_k$ we deduce that
$Q_{\mfp_a\mfn_k}(X',Y')=\Phi(C_{\mfp_a}X',Y')$. From Lemma
\ref{traces1}, we also obtain
$Q_{\mfp_a\mfn_j}(X',Y')=Q_{\mfn_j\mfp_a}(X',Y')=0$, for $j\neq k$ and
$Q_{\mfn_k\mfp_a}(X',Y')=Q_{\mfp_a\mfn_k}(X',Y')=\Phi(C_{\mfp_a}X',Y')$,
for $j=k$.

$\Box$ \eproof

In the remaining of this section, we obtain formulae for the Ricci curvature of an adapted metric of the fibration $M=G/L\rightarrow G/K=N$ in the vertical, horizontal and $\mfp\times\mfn$ directions, by using Lemma \ref{tracespn} and the formula for the Ricci curvature from Theorem \ref{ric1}.

\subsection{The Ricci Curvature of an Adapted metric in the Vertical Direction}

\blem\label{ricF} Let $g_F$ be the $K$-invariant metric on the fiber $F=K/L$ determined by the $Ad\,L$-invariant Euclidean product $g_{\mfp}=\oplus_{a=1}^s\la_aB_{\mfp_a}$. The Ricci curvature of $g_F$ is given by
$Ric^F=\oplus_{a=1}^sq_aB_{\mfp_a}$, where

\beq\label{constant8}q_a=\frac{1}{2}\sum_{b,c=1}^s\left(\frac{\la_a^2}{2\la_c\la_b}-\frac{\la_c}{\la_b}\right)q^{cb}_a+\frac{\ga_a}{2}.\eeq

The constants $q^{cb}_a$ and $\ga_a$ are such that
\beqar\label{constant1} \Phi_{\mfk}\mid_{\mfp_a\times\mfp_a}=\ga_a\Phi\mid_{\mfp_a\times\mfp_a}\\
\label{constant2} Q_{\mfp_b\mfp_c}\mid_{\mfp_a\times\mfp_a}=q^{cb}_a\Phi\mid_{\mfp_a\times\mfp_a}.\eeqar

In particular, $Ric^F(\mfp_a,\mfp_b)=0$, if $a\neq b$. \elem

\bproof  Since $\mfp_1,\ldots,\mfp_s$ are pairwise inequivalent
irreducible $Ad\,L$-submodules and the Ricci curvature of $g_F$, $Ric^F$, is an $Ad\,L$-invariant
symmetric bilinear form, we may write
$Ric^F=\oplus_{a=1}^sq_aB_{\mfp_a}$, for some
constants $q_1,\ldots,q_s$. In particular, we have
$Ric^F(\mfp_a,\mfp_b)=0$, if $a\neq b$. By Theorem \ref{ric1}, the Ricci curvature of $g_F$ is

$$Ric^F(X,X)=\dfrac{1}{2}\sum_{b,c=1}^s\left(\frac{\la_b}{\la_c}-\dfrac{\la_a^2}{2\la_c\la_b}\right)Q_{\mfp_c\mfp_b}(X,X)-\frac{1}{2}\Phi_{\mfk}(X,X).$$

By Lemma \ref{traces1} the maps $Q_{\mfp_c\mfp_b}$ are
$Ad\,L$-invariant symmetric bilinear maps. Since $\mfp_a$ is
$Ad\,L$-irreducible, there are constants $q^{cb}_a$ as defined by (\ref{constant2}). Similarly, there is a constant $\ga_a$ as defined in (\ref{constant1}). By the expression above for $Ric^F$, we must have

$$q_a=\frac{1}{2}\sum_{b,c=1}^s\left(\frac{\la_a^2}{2\la_c\la_b}-\frac{\la_c}{\la_b}\right)q^{cb}_a+\frac{\ga_a}{2}$$

and the result follows from this.

$\Box$ \eproof

\li

\bprop \label{riccip}Let $g_M$ be an adapted metric on the homogeneous fibration $M=G/L\rightarrow G/K$ defined by the Euclidean product

$$g_{\mfm}=\left(\oplus_{a=1}^s\la_aB_{\mfp_a}\right)\oplus\left(\oplus_{k=1}^n\mu_kB_{\mfn_k}\right).$$

We have $Ric(\mfp_a,\mfp_b)=0$, if $a\neq b$. For $X\in\mfp_a$,

$$Ric(X,X)=\left(q_a+\dfrac{\la_a^2}{4}\sum_{j=1}^n\dfrac{c_{\mfn_j,a}}{\mu_j^2}\right)B(X,X),$$

where, for $j=1,\ldots,n$, the constants $c_{\mfn_j,a}$ are such
that

\beq\label{constant3}\Phi(C_{\mfn_j}\cdot,\cdot)\mid_{\mfp_a\times\mfp_a}=c_{\mfn_j,a}\Phi\mid_{\mfp_a\times\mfp_a}\eeq

and $C_{\mfn_j}$ is the Casimir operator of $\mfn_j$ with respect
to $\Phi$. The constant $q_a$ is such that

$$q_a=\frac{1}{2}\sum_{b,c=1}^s\left(\frac{\la_a^2}{2\la_c\la_b}-\frac{\la_c}{\la_b}\right)q^{cb}_a+\frac{\ga_a}{2},$$

with $q^{cb}_a$ and $\ga_a$ defined by

\beqar\label{constant4}\Phi_{\mfk}\mid_{\mfp_a\times\mfp_a}=\ga_a\Phi\mid_{\mfp_a\times\mfp_a}\\
\label{constant5}Q_{\mfp_b\mfp_c}\mid_{\mfp_a\times\mfp_a}=q^{cb}_a\Phi\mid_{\mfp_a\times\mfp_a}.\eeqar\eprop

\bproof Since $\mfp_1,\ldots,\mfp_s$ are pairwise inequivalent
irreducible $Ad\,L$-submodules and $Ric\mid_{\mfp\times\mfp}$ is
an $Ad\,L$-invariant symmetric bilinear form, we have that
$Ric\mid_{\mfp\times\mfp}$ is diagonal, i.e.,

$$Ric\mid_{\mfp\times\mfp}=a_1B_{\mfp_1}\oplus\ldots\oplus a_sB_{\mfp_s},$$

for some constants $a_1,\ldots,a_s$. In particular, we have
$Ric(\mfp_a,\mfp_b)=0$, if $a\neq b$. Hence, $Ric\mid_{\mfp\times\mfp}$ is determined by elements
$Ric(X,X)$ with $X\in\mfp_a$, $a=1,\ldots,s$.

\li

By Lemma \ref{tracespn} we obtain that only
$Q_{\mfn_j\mfn_j}(X,X)=\Phi(C_{\mfn_j}X,X)$ and
$Q_{\mfp_b\mfp_c}(X,X)$ is non-zero. Therefore, by Theorem
\ref{ric1} we obtain that

$$Ric(X,X)=$$

$$\frac{1}{2}\sum_{j,k=1}^s\left(\frac{\la_k}{\la_j}-\frac{\la_a^2}{2\la_j\la_k}\right)Q_{\mfp_j\mfp_k}(X,X)+\frac{1}{2}\sum_{j=1}^n\left(1-\frac{\la_a^2}{2\mu_j^2}\right)Q_{\mfn_j\mfn_j}(X,X)-\frac{1}{2}\Phi(X,X).$$

\li

We have $\sum_{j=1}^mQ_{\mfn_j\mfn_j}(X,X) = \sum_{j=1}^m\Phi(C_{\mfn_j}X,Y) = \Phi(X,X)-\Phi_{\mfk}(X,X)$.

Hence we can rewrite $Ric(X,X)$ as follows:

$$\underbrace{\frac{1}{2}\sum_{a,b=1}^s\left(\frac{\la_b}{\la_c}-\frac{\la_a^2}{2\la_c\la_b}\right)Q_{\mfp_c\mfp_b}(X,X)-\frac{1}{2}\Phi_{\mfk}(X,X)}_{(1)}-\frac{1}{2}\sum_{j=1}^n\frac{\la_a^2}{2\mu_j^2}\Phi(C_{\mfn_j}X,X).$$

As we saw in the proof of Lemma \ref{ricF}, the summand (1) is
just $Ric^F(X,X)=q_aB(X,X)$. Furthermore, since $\Phi(C_{\mfn_j}\cdot,\cdot)=Q_{\mfn_j\mfn_j}$ there are constants
$c_{\mfn_j,a}$, defined in (\ref{constant3}).

Therefore,

$$Ric(X,X)=\left(q_a+\frac{1}{4}\sum_{j=1}^n\frac{\la_a^2}{\mu_j^2}c_{\mfn_j,a}\right)B(X,X).$$

$\Box$\eproof

\li

\subsection{The Ricci Curvature of an Adapted metric in the Horizontal Direction}

\blem\label{ricN}Let $g_N$ be the $G$-invariant metric on $N=G/K$ determined by the $Ad\,L$-invariant Euclidean product $g_{\mfn}=\oplus_{j=1}^n\mu_jB_{\mfn_j}$. The Ricci curvature of $g_N$ is given by $Ric^N=\oplus_{k=1}^nr_kB_{\mfn_k}$, where

\beq\label{constant7}r_k=\frac{1}{2}\sum_{j,i=1}^n\left(\frac{\mu_a^2}{2\mu_i\mu_j}-\frac{\mu_i}{\mu_j}\right)r^{ji}_k+\frac{1}{2}\eeq

and the constants $r^{ji}_k$ are such that

\beq\label{constant8}Q_{\mfn_j\mfn_i}\mid_{\mfn_k\times\mfn_k}=r^{ji}_k\Phi\mid_{\mfn_k\times\mfn_k}.\eeq

In particular, $Ric^N(\mfn_k,\mfn_j)=0$, if $k\neq j$. \elem

\bproof The proof is similar to the proof of Lemma \ref{ricF} and it can be found in detail in \cite{Fa}.

$\Box$\eproof

\li

\bprop \label{riccin} Let $g_M$ be an adapted metric on the homogeneous fibration $M=G/L\rightarrow G/K$ defined by the Euclidean product

$$g_{\mfm}=\left(\oplus_{a=1}^s\la_aB_{\mfp_a}\right)\oplus\left(\oplus_{k=1}^n\mu_kB_{\mfn_k}\right).$$

We have $Ric(\mfn_k,\mfn_j)=0$, if $j\neq k$. For $X\in\mfn_k$,

$$Ric(X,X)=-\dfrac{1}{2\mu_k}\sum_{a=1}^s\la_aB(C_{\mfp_a}X,X)+r_kB(X,X),$$

where $C_{\mfp_a}$ is the Casimir
operator of $\mfp_a$ with respect to $\Phi$,

$$r_k=\frac{1}{2}\sum_{j,i=1}^n\left(\frac{\mu_a^2}{2\mu_i\mu_j}-\frac{\mu_i}{\mu_j}\right)r^{ji}_k+\frac{1}{2}$$

and the constants $r^{ji}_k$ are such that

$$Q_{\mfn_j\mfn_i}\mid_{\mfn_k\times\mfn_k}=r^{ji}_k\Phi\mid_{\mfn_k\times\mfn_k}.$$\eprop

\bproof Let $X\in\mfn_k$ and $Y\in\mfn_{k'}$. By Lemma
\ref{tracespn} we have $Q_{\mfp_a\mfp_b}(X,Y)=0$, for every
$a,b=1,\ldots,s$. Also,
$Q_{\mfn_j\mfp_a}(X,Y)=Q_{\mfp_a\mfn_j}(X,Y)=0$, if $j\neq k,k'$.
Therefore, it follows from Theorem \ref{ric1} and Lemma \ref{ricN}
that, if $k\neq k'$, then

$$Ric(X,Y)=\frac{1}{2}\sum_{j,i=1}^n\left(\frac{\mu_i}{\mu_j}-\frac{\mu_k\mu_{k'}}{2\mu_i\mu_j}\right)Q_{\mfn_j\mfn_i}(X,Y)-\frac{1}{2}\Phi(X,Y)=Ric^N(X,Y)=0.$$

Hence, $Ric\mid_{\mfn\times\mfn}$ is determined by elements
$Ric(X,X)$ with $X\in\mfn_k$, $k=1,\ldots,n$. For $X\in\mfn_k$, by
Theorem \ref{ric1}, we get

$$Ric(X,X)=$$

$$\frac{1}{2}\sum_{k=1}^s\left(\frac{\mu_k}{\la_a}-\frac{\mu_k^2}{2\mu_
k\la_a}\right)Q_{\mfp_a\mfn_k}(X,X)+\frac{1}{2}\sum_{a=1}^s\left(\frac{\la_a}{\mu_k}-\frac{\mu_k^2}{2\mu_
k\la_a}\right)Q_{\mfn_k\mfp_a}(X,X)+Ric^N(X,X).$$

From Lemma \ref{tracespn}, we know that
$Q_{\mfn_k\mfp_a}(X,X)=Q_{\mfp_a\mfn_k}(X,X)=\Phi(C_{\mfp_a}X,X)$.
Hence, we simplify the expression above obtaining

$$Ric(X,X)=\dfrac{1}{2}\sum_{a=1}^s\dfrac{\la_a}{\mu_k}\Phi(C_{\mfp_a}X,X)+Ric^N(X,X).$$
Finally, by using Lemma \ref{ricN} we have $Ric^N(X,X)=r_kB(X,X)$ and this concludes the proof.

$\Box$\eproof

\li

\subsection{The Ricci Curvature of an Adapted metric in the Mixed Direction $\mfp\times\mfn$}

\bprop \label{riccipn} Let $g_M$ be an adapted metric on the homogeneous fibration $M=G/L\rightarrow G/K$ defined by the Euclidean product

$$g_{\mfm}=\left(\oplus_{a=1}^s\la_aB_{\mfp_a}\right)\oplus\left(\oplus_{k=1}^n\mu_kB_{\mfn_k}\right).$$

For $X\in\mfp_a$, $Y\in\mfn_k$,
$$Ric(X,Y)=\dfrac{\la_a\mu_k}{4}\sum_{j=1}^n\dfrac{B(C_{\mfn_j}X,Y)}{\mu_j^2},$$

where $C_{\mfn_j}$ is the Casimir
operator of $\mfn_j$ with respect to $\Phi$. \eprop

\bproof For $X\in\mfp$ we know from Lemma \ref{tracespn} that
$Q_{\mfn_j\mfp_a}(X,Y)=Q_{\mfp_a\mfn_j}(X,Y)=0$ and $Q_{\mfn_i\mfn_j}(X,Y)=0$, if
$i\neq j$, whereas
$Q_{\mfn_j\mfn_j}(X,Y)=\Phi(C_{\mfn_j}X,Y)$.
Moreover, for $Y\in\mfn_k$, since
$ad_Xad_Y\mfp\subset\mfn_k\perp\mfp$, we
also obtain from Lemma \ref{tracespn} that $Q_{\mfp_b\mfp_c}(X,Y)=0$. Therefore, only
$Q_{\mfn_j\mfn_j}(X,Y)=\Phi(C_{\mfn_j}X,Y)$, may
not be zero. Furthermore, $\Phi(X,Y)=0$. Hence, from Theorem
\ref{ric1} we get

$$Ric(X,Y)=\frac{1}{2}\sum_{j=1}^n\left(1-\frac{\la_a\mu_k}{2\mu_j^2}\right)Q_{\mfn_j\mfn_j}(X,Y)=\frac{1}{2}\sum_{j=1}^n\left(1-\frac{\la_a\mu_k}{2\mu_j^2}\right)\Phi(C_{\mfn_j}X,Y).$$

On the other hand,
$$\sum_{j=1}^n\Phi(C_{\mfn_j}X,Y)=\Phi(C_{\mfn}X,Y)=\Phi(X,Y)-\Phi(C_{\mfk}X,Y)=0,$$
since $C_{\mfk}\mfp\subset\mfk\perp\mfn$. Therefore,

$$Ric(X,Y)=-\frac{\la_a\mu_k}{4}\sum_{j=1}^n\frac{\Phi(C_{\mfn_j}X,Y)}{\mu_j^2}.$$

$\Box$\eproof

\li

\subsection{Necessary Conditions for the Existence of an Adapted Einstein Metric}

From the expressions obtained previously for the Ricci curvature
in the horizontal direction and in the direction of
$\mfp\times\mfn$ we obtain two necessary conditions for the
existence of an adapted Einstein metric on $M$.

\bcor \label{cond1}Let $g_M$ be an adapted metric on the homogeneous fibration $M=G/L\rightarrow G/K$ defined by the Euclidean product

$$g_{\mfm}=\left(\oplus_{a=1}^s\la_aB_{\mfp_a}\right)\oplus\left(\oplus_{k=1}^n\mu_kB_{\mfn_k}\right).$$

If $g_M$ is Einstein, then the operator $\sum_{a=1}^s\la_aC_{\mfp_a}$ is scalar on each $\mfn_j$. \ecor

\bproof Let $g_M$ be an adapted metric as defined in (\ref{mdef}).
If $g_M$ is Einstein with Einstein constant $E$, then,
$Ric\mid_{\mfn\times\mfn}=Eg_{\mfn}$ and thus
$Ric\mid_{\mfn\times\mfn}$ is $Ad\,K$-invariant. Therefore, by
Proposition \ref{riccin}, we conclude that
$\sum_{a=1}^s\la_aC_{\mfp_a}\mid_{\mfn}$ must be
$Ad\,K$-invariant. Hence, $\sum_{a=1}^s\la_aC_{\mfp_a}\mid_{\mfn_k}$
is scalar, by irreducibility of $\mfn_k$.

$\Box$ \eproof

\bcor \label{cond2}Let $g_M$ be an adapted metric on the homogeneous fibration $M=G/L\rightarrow G/K$ defined by the Euclidean product

$$g_{\mfm}=\left(\oplus_{a=1}^s\la_aB_{\mfp_a}\right)\oplus\left(\oplus_{k=1}^n\mu_kB_{\mfn_k}\right).$$
The orthogonality condition $Ric(\mfp,\mfn)=0$
holds if and only if

\beq\label{ort1}\sum_{j=1}^n\frac{1}{\mu_j^2}C_{\mfn_j}(\mfp)\subset
\mfk.\eeq

Moreover, if $g_M$ is Einstein, then (\ref{ort1}) holds. \ecor

\bproof From Proposition
\ref{riccipn}, we obtain that $Ric(\mfp,\mfn)=0$ if and only if,
for every $X\in\mfp_a$ and $Y\in\mfn_b$,

$$\Phi\left(\sum_{j=1}^n\frac{C_{\mfn_j}}{\mu_j^2}X,Y\right)=0.$$

This holds if only if
$\sum_{j=1}^n\frac{C_{\mfn_j}}{\mu_j^2}X\subset \mfk$, for every
$X\in\mfp$.

If $g_M$ is Einstein with Einstein constant $E$, then
$Ric(\mfp,\mfn)=Eg_{\mfm}(\mfp,\mfn)=0$.

$\Box$ \eproof

The two previous Corollaries may be restated as in Theorem \ref{cond3},
which emphasizes the fact that the two necessary conditions
obtained for existence of an Einstein adapted metric are just
algebraic conditions on the Casimir operators of the submodules
$\mfp_a$ and $\mfn_k$.

\li

\textbf{\emph{Proof of Theorem \ref{cond3}:}} The assertions follow from Corollaries \ref{cond1}
 and \ref{cond2}.  In \ref{cond2} we set $\nu_k=1/\mu_k^2$. Hence, $\nu_1=\ldots=\nu_n$ occurs when $g_N$ is the standard metric. Moreover, if $\nu_1=\ldots=\nu_n$, the inclusion in Theorem \ref{cond2} is equivalent to $C_{\mfn}(\mfp)\subset\mfk$, which always holds since $C_{\mfn}=Id-C_{\mfk}$ and $C_{\mfk}$ maps $\mfp$ into $\mfk$. So we obtain a condition on the $C_{\mfn_j}$'s only when $g_N$ is not standard.

 $\Box$

\section{Einstein Binormal Metrics}\label{sectionBinormal}

In this Section we obtain the Ricci curvature of a binormal metric ( see (\ref{mdefbi})) and conditions for such a metric to be Einstein. We prove Theorem \ref{binormal1} and the subsequent Corollaries stated in Section \ref{intro}. As we shall see, the conditions for the existence of an Einstein binormal metric translate in very simple conditions on the Casimir
operators of $\mfk$, $\mfl$ and $\mfp_a$. As in previous sections, we consider the homogeneous fibration $M=G/L\rightarrow G/K=N$ as in (\ref{fibDefIntro}). We use the notation from Sections \ref{intro} and \ref{sectionRF}. We start by considering the case when the restriction $g_F$ of the metric to the fiber $F=K/L$ is normal.

\bprop \label{riccinf} Let $g_M$ be an adapted metric on the homogeneous fibration $M=G/L\rightarrow G/K$ defined by the Euclidean product

$$g_{\mfm}=\la B_{\mfp}\oplus\left(\oplus_{k=1}^n\mu_kB_{\mfn_k}\right).$$

The Ricci curvature of $g_M$ is given by:

(i) For $X\in\mfp_a$,

$$Ric(X,X)=\left(q_a+\frac{\la^2}{4}\sum_{j=1}^n\dfrac{c_{\mfn_j,a}}{\mu_j^2}\right)B(X,X),$$

where $q_a=\frac{1}{2}\left(c_{\mfl,a}+\frac{\ga_a}{2}\right)$, $c_{\mfl,a}$ is the eigenvalue of the Casimir operator of
$\mfl$ on $\mfp_a$, $\ga_a$ and $c_{\mfn_j,a}$ are given by (\ref{constant1}) and (\ref{constant3}), respectively.

(ii) For $X\in\mfn_k$,

$$Ric(X,X)=-\frac{\la}{2\mu_k}B(C_{\mfp}X,X)+r_kB(X,X),$$

where $r_k$ is given by (\ref{constant7});

(iii) For $X\in\mfp_a$ and $Y\in\mfn_k$,
$$Ric(X,Y)=\dfrac{\la\mu_k}{4}\sum_{j=1}^n\dfrac{B(C_{\mfn_j}X,Y)}{\mu_j^2};$$

(iv) $Ric(\mfp_a,\mfp_b)=0$, if $a\neq b$, and $Ric(\mfn_i,\mfn_j)=0$, if $i\neq j$. \eprop

\bproof (i) By Corollary \ref{ric1cor1}, if $g_F$ is a multiple of
$B=-\Phi$, then we obtain that

$$Ric^F(X,X)=-\frac{1}{4}\Phi_{\mfk}(X,Y)-\frac{1}{2}\Phi(C_{\mfl}X,Y)= -\frac{1}{2}\left(\frac{\ga_a}{2}+c_{\mfl,a}\right)\Phi(X,X),$$

for $X\in\mfp_a$, where $c_{\mfl,a}$ is the eigenvalue of the
Casimir operator of $\mfl$ with respect to $\Phi$ on $\mfp_a$. Hence, we have

$$q_a=\frac{1}{2}\left(\frac{\ga_a}{2}+c_{\mfl,a}\right).$$

The result then follows from this and Proposition \ref{riccip}.

\li

(ii) It follows directly from Proposition \ref{riccin}, by observing
that $\sum_{a=1}^sC_{\mfp_a}=C_{\mfp}$.

(iii) The Ricci curvature in the direction $\mfp\times\mfn$
essentially remains unchanged; the expression given is just that of
Proposition \ref{riccipn} after replacing $\la_1,\ldots,\la_s$ by
$\la$.

(iv) These orthogonality conditions are satisfied by any adapted
metric on $M$ and were shown to hold in Propositions \ref{riccip}
and \ref{riccin}.

$\Box$ \eproof

Similarly, if the metric on the base is normal, then we obtain the following characterization:

\bprop\label{riccinb} Let $g_M$ be an adapted metric on the homogeneous fibration $M=G/L\rightarrow G/K$ defined by the Euclidean product
$$g_{\mfm}=\left(\oplus_{a=1}^s\la_aB_{\mfp_a}\right)\oplus\mu B_{\mfn}.$$
The Ricci curvature of $g_M$ is given by:

(i) For $X\in\mfp_a$,

$$Ric(X,X)=\left(q_a+\frac{\la_a^2}{4\mu^2}(1-\ga_a)\right)B(X,X),$$

where $q_a$ and $\ga_a$ are given by (\ref{constant1}) and (\ref{constant8}), respectively.

(ii) For $X\in\mfn_k$,
$$Ric(X,X)=-\dfrac{1}{2}\sum_{a=1}^s\dfrac{\la_a}{\mu}B(C_{\mfp_a}X,X)+r_kB(X,X),$$

with $r_k=\frac{1}{2}\left(\frac{1}{2}+c_{\mfk,k}\right)$, where $c_{\mfk,k}$ is the eigenvalue of the Casimir operator
$C_{\mfk}$ on $\mfn_k$;

(iii) $Ric(\mfp,\mfn)=0$;

(iv) $Ric(\mfp_a,\mfp_b)=0$, if $a\neq b$, and $Ric(\mfn_i,\mfn_j)=0$, if $i\neq j$. \eprop

\bproof (i) From the fact that $\ga_a+\sum_{j=1}^n c_{\mfn_j,a}=1$,
we obtain

$$\sum_{j=1}^n\frac{\la_a^2}{\mu^2}c_{\mfn_j,a}=\frac{\la_a^2}{\mu^2}(1-\ga_a).$$

The required expression follows immediately from Proposition
\ref{riccip}.

\li

(ii) From Corollary \ref{ric1cor1} we obtain that
$r_k=\frac{1}{2}\left(\frac{1}{2}+c_{\mfk,k}\right)$, where $r_k$ is
as defined in Lemma \ref{ricN}. The expression then follows from
Proposition \ref{riccin}.

\li

(iii) By using the fact that $C_{\mfn}=\sum_{j=1}^nC_{\mfn_j}$, from
Proposition \ref{riccipn} it follows that

$$Ric(X,Y)=\frac{\la_a}{4\mu}B(C_{\mfn}X,Y),$$

for every $X\in\mfp_a$ and $Y\in\mfn_k$. Moreover, since
$C_{\mfn}=C_{\mfg}-C_{\mfk}=Id-C_{\mfk}$ and
$C_{\mfk}(\mfp)\subset\mfk$, we have that $C_{\mfn}(X)\in\mfk$ is
orthogonal to $Y\in\mfn$ with respect to $B$. Hence, $Ric(X,Y)=0$.

\li

(iv) these orthogonality conditions are simply those in
Propositions \ref{riccip} and \ref{riccin}.

$\Box$\eproof

The Ricci curvature of a binormal metric follows immediately from Propositions \ref{riccinf} and \ref{riccinb}:

\bcor \label{riccitypeI} Consider the homogeneous fibration $M=G/L\rightarrow G/K$ and the binormal metric $g_M$ on $M$ defined by the Euclidean product $g_{\mfm}=\la B_{\mfp}\oplus\mu B_{\mfn}$.

(i) For every $X\in\mfp_a$,

$$Ric(X,X)=\left(q_a+\frac{\la^2}{4\mu^2}(1-\ga_a)\right)B(X,X),$$

where $q_a=\frac{1}{2}\left(\frac{\ga_a}{2}+c_{\mfl,a}\right)$,
$c_{\mfl,a}$ is the eigenvalue of $C_{\mfl}$ on $\mfp_a$ and $\ga_a$
is determined by

$$\Phi_{\mfk}\mid_{\mfp_a\times\mfp_a}=\ga_a\Phi\mid_{\mfp_a\times\mfp_a};$$

(ii) For every $Y\in\mfn_j$,
$$Ric(Y,Y)=-\frac{\la}{2\mu}B(C_{\mfp}Y,Y)+r_jB(Y,Y),$$

where $r_j=\frac{1}{2}\left(\frac{1}{2}+c_{\mfk,j}\right)$ and
$c_{\mfk,j}$ is the eigenvalue of $C_{\mfk}$ on $\mfn_j$;

(iii) $Ric(\mfp,\mfn)=0$;

(iv) $Ric(\mfp_a,\mfp_b)=0$, if $a\neq b$, and $Ric(\mfn_i,\mfn_j)=0$, if $i\neq j$.

\ecor

We finally prove Theorem \ref{binormal1} given in Section \ref{intro}.

\li

\textbf{\emph{Proof of Theorem \ref{binormal1}:}} Let $g_M$ be an adapted metric on $M$ and $g_{\mfm}=\left(\oplus_{a=1}^s\la_aB_{\mfp_a}\right)\oplus\left(\oplus_{k=1}^n\mu_kB_{\mfn_k}\right)$ the associated $Ad\,L$-invariant Euclidean product on $\mfm$. By Corollary \ref{cond1}, we have that, if $g_M$ is Einstein, then
$C_{\mfp}$ and $C_{\mfl}$ are scalar on $\mfn_j$, for every
$j=1,\ldots,n$. Say

$$C_{\mfp}\mid_{\mfn_j}=b^jId \textrm{ and } C_{\mfl}\mid_{\mfn_j}=c_{\mfl,j}Id.$$

Suppose that $g$ is Einstein with constant $E$. From Corollary
\ref{riccitypeI}, we obtain the Einstein equations

\beqar \label{e1}-\dfrac{\la}{2\mu}b^j+r_j=\mu E,\,j=1,\ldots,n\\
\label{e2}\dfrac{1}{2}\left(\frac{\ga_a}{2}+c_{\mfl,a}+\frac{\la^2}{2\mu^2}(1-\ga_a)\right)=\la
E,\,a=1,\ldots,s.\eeqar

\li

If $n>1$, from Equation (\ref{e1}) we obtain the following:

\beq\label{difij}
\dfrac{\la}{2\mu}(b^i-b^j)=r_i-r_j,\,i,j=1,\ldots,n.\eeq

By using Lemma \ref{ricN} we have
$$r_i-r_j=\dfrac{1}{2}\left(\dfrac{1}{2}+c_{\mfk,i}\right)-\dfrac{1}{2}\left(\dfrac{1}{2}+c_{\mfk,j}\right)=\dfrac{1}{2}(c_{\mfk,i}-c_{\mfk,j}),$$

whereas

$$b^i-b^j=(c_{\mfk,i}-c_{\mfk,j})-(c_{\mfl,i}-c_{\mfl,j}).$$

Therefore, Equation (\ref{difij}) becomes
$$-\frac{\la}{\mu}(\underbrace{c_{\mfl,i}-c_{\mfl,j}}_{\de_{ij}^{\mfl}})=\left(1-\frac{\la}{\mu}\right)(\underbrace{c_{\mfk,i}-c_{\mfk,j}}_{\de_{ij}^{\mfk}}).$$

By using the variable $X$, we rewrite the equation above as
$-\frac{1}{X}\de_{ij}^{\mfl}=\left(1-\frac{1}{X}\right)\de_{ij}^{\mfk}$,
and this yields $\de_{ij}^{\mfl}=(1-X)\de_{ij}^{\mfk}$.

\li

Equation (\ref{e2}) may be rewritten as

\beq\label{e3}\dfrac{1}{2}\left(\frac{\ga_a}{2}+c_{\mfl,a}\right)X+(1-\ga_a)\frac{1}{4X}=\mu
E.\eeq

Hence, if $s>1$, for $a,b=1,\ldots,s$, we get

$$\dfrac{1}{2}\left(\frac{\ga_a}{2}+c_{\mfl,a}\right)X+(1-\ga_a)\frac{1}{4X}=\dfrac{1}{2}\left(\frac{\ga_b}{2}+c_{\mfl,b}\right)X+(1-\ga_b)\frac{1}{4X},$$

which yields

$$\underbrace{c_{\mfl,a}-c_{\mfl,b}}_{\de_{ab}^{\mfl}}=\frac{1}{2}\left(\frac{1}{X^2}-1\right)(\underbrace{\ga_a-\ga_b}_{\de_{ab}^{\mfk}}).$$

By solving this equation we obtain

$$(2\de_{ab}^{\mfl}+\de_{ab}^{\mfk})X^2=\de_{ab}^{\mfk}.$$

\li

Finally, by using Equations (\ref{e1}) and (\ref{e3}) we obtain the
equality

$$\dfrac{1}{2}\left(\frac{\ga_a}{2}+c_{\mfl,a}\right)X+(1-\ga_a)\frac{1}{4X}=-\dfrac{b^j}{2X}+\frac{1}{2}\left(\frac{1}{2}+c_{\mfk,j}\right),$$

which rearranged gives

$$\left(\frac{\ga_a}{2}+c_{\mfl,a}\right)X^2-\left(\frac{1}{2}+c_{\mfk,j}\right)X+\frac{1}{2}(1-\ga_a+2b^j)=0.$$

$\Box$

\li

\textbf{\emph{Proof of Corollary \ref{bincor1}:}} Since $\mfp$ is an irreducible $Ad\,L$-module and $\mfn$
is an irreducible $Ad\,K$-module, then any adapted metric on $M$
is binormal. Hence, we use Theorem \ref{binormal1}. By the
irreducibility of $\mfp$ and $\mfn$, we have $s=1$ and $n=1$ and
thus Einstein binormal metrics are given by positive solutions of
(\ref{einI3}), if $C_{\mfp}$ is scalar on $\mfn$. Hence, from
Theorem \ref{binormal1} we conclude that there exists on $M$ an
Einstein binormal metric if and only if $C_{\mfp}$ is scalar on
$\mfn$ and $\triangle\geq 0$, where
$$\triangle=(1+2c_{\mfk,\mfn})^2-4(\ga+2c_{\mfl,\mfp})(1-\ga+2b).$$

\li

Since $F$ is isotropy irreducible and $dim\,F>1$, we have
$\ga+2c_{\mfl,p}\neq 0$ and the polynomial in (\ref{einI3}) has
exactly degree two. In fact, if $\ga+2c_{\mfl,p}=0$, then
$\ga=c_{\mfl,p}=0$ and thus, in particular, $\mfp$ lies in the
center of $\mfk$. But the hypothesis that $\mfp$ is irreducible
and abelian implies that $\mfp$ is $1$-dimensional which
contradicts the hypothesis that $dim\,F>1$. Therefore,
$\ga+2c_{\mfl,p}\neq 0$. In this case, the solutions of
(\ref{einI3}) are

$$X=\frac{1+2c_{\mfk,\mfn}\pm\sqrt{\triangle}}{2(\ga+2c_{\mfl,\mfp})}.$$

$\Box$

\li

\emph{\textbf{Proof of Corollary \ref{bincor2}:}} The fact that $\mfp$ is 1-dimensional implies that $\mfp$
lies in the center of $\mfk$. Hence, in the notation of Corollary
\ref{bincor1}, $\ga=c_{\mfl,\mfp}=0$. On the other hand, if $\mfn$
is $Ad\,K$-irreducible then, the semisimple part of $K$ acts
transitively on $\mfn$. Moreover, since $\mfp$ lies in the center
of $\mfk$, then the semisimple part of $\mfl$ coincides with the
semisimple part of $\mfk$. Hence, $L$ also acts transitively on
$\mfn$ and $\mfn$ is an irreducible $Ad\,L$-module as well.
Consequently, any $G$-invariant metric on $M$ is adapted and
moreover is binormal, by the irreducibility of $\mfp$ and $\mfn$.
Furthermore, $C_{\mfp}$ must be scalar on $\mfn$, since $C_{\mfk}$
and $C_{\mfl}$ are scalar on $\mfn$. Therefore, $G$-invariant
Einstein metrics are given by positive solutions of (\ref{einI3})
in Theorem \ref{binormal1}. Since $\ga=c_{\mfl,\mfp}=0$,
(\ref{einI3}) is just a degree-one equation whose solution is

\beq\label{solcircle1} X=\frac{1+2b}{1+2c_{\mfk,\mfn}},\eeq

where $b$ is the eigenvalue of $C_{\mfp}$ on $\mfn$. Since $\mfg$ is simple we have $tr(C_{\mfp})=dim\,\mfp=1$. Since
$\mfp$ lies in the center of $\mfk$, $C_{\mfp}$ vanishes on $\mfk$
and thus
$tr(C_{\mfp})=tr(C_{\mfp}\mid_{\mfn})=b\,\textrm{dim}\,\mfn=bm$.
Hence, $b=1/m$. By replacing $b$ on (\ref{solcircle1}) we obtain
the expression given for $X$.

$\Box$

\li

\textbf{\emph{Proof of Corollary \ref{binormal3}:}} If $\Phi\circ C_{\mfl}\mid_{\mfp\times\mfp}=\al
\Phi_{\mfk}\mid_{\mfp\times\mfp}$, then
\beq\label{eqcas1}c_{\mfl,a}=\al \ga_a, \textrm{ for every
}a=1,\ldots,s.\eeq

Therefore, for any $a,b=1,\ldots,s$, if $s>1$,
$2\de_{ab}^{\mfl}+\de_{ab}^{\mfk}=(2\al+1)\de_{ab}^{\mfk}$ and,
thus, Equation (\ref{einI2}) in Theorem \ref{binormal1} becomes

\beq\label{eqcas2}(2\al+1)\de_{ab}^{\mfk}X^2=\de_{ab}^{\mfk}.\eeq

In particular, (\ref{eqcas1}) implies that $c_{\mfl,a}=0$ if and
only if $\ga_a=0$ (thus if $\mfp$ has submodules where $L$ acts
trivially, then $\Phi_{\mfk}$ vanish on those submodules and then
they lie in the center of $\mfk$. If $K$ is semisimple, then the
isotropy representation of $K/L$ is faithful). The fact that the
isotropy representation of $K/L$ is not irreducible implies that
$\mfp$ decomposes as a direct sum $\mfp_1\oplus\ldots\oplus\mfp_s$
with $s>1$. For the indices for which
$\ga_a\neq\ga_b$, we have $\de_{ab}^{\mfk}\neq 0$ and (\ref{eqcas2})
implies that
$$X=\frac{1}{\sqrt{2\al+1}}.$$

Hence, $X=\frac{1}{\sqrt{2\al+1}}$ must be a root of the polynomial
in (\ref{einI3}). By using the fact that $c_{\mfl,a}=\al \ga_a$ and
$b_j=c_{\mfk,j}-c_{\mfl,j}$, simple calculations show that

\beq\label{casj1}c_{\mfl,j}=\left(1-\frac{1}{\sqrt{2\al+1}}\right)\left(c_{\mfk,j}+\frac{1}{2}
\right).\eeq

We observe that this condition implies (\ref{einI1}) in Theorem
\ref{binormal1}, as we can see by the equalities below:

$$\de_{ij}^{\mfl}=c_{\mfl,i}-c_{\mfl,j}=\left(1-\frac{1}{\sqrt{2\al+1}}\right)(c_{\mfk,i}-c_{\mfk,j})=(1-X)\de_{ij}^{\mfk}.$$

Hence, there is a binormal Einstein metric if and only if
(\ref{casj1}) is satisfied and the operator $C_{\mfp}$ is scalar
on $\mfn_j$, for every $j=1,\ldots,n$. In this case, according
also to Theorem \ref{binormal1} such metric is, up to homothety,
given by $B_{\mfp}\oplus
\frac{1}{\sqrt{2\al+1}}B_{\mfn}$.

\li

It remains to show that if exists an Einstein binormal metric then $\sqrt{2\al+1}\in\rationals$. This follows from the fact that the eigenvalues of $C_{\mfk}$ and $C_{\mfl}$ on $\mfn_j$ are rational numbers. Since
$\mfk$ is a compact algebra, the eigenvalue of its Casimir operator
on the complex representation on $\mfn_j^{\complex}$ is given by

$$\frac{<\la_j,\la_j+2\de>}{2h^*(\mfg)}\in \rationals,$$

where $\la_j$ is the highest weight for  $\mfn_j^{\complex}$,
$2\de$ is the sum of all positive roots of $\mfk$ and $h^*(\mfg)$
is the dual Coxeter number of $\mfg$ (\cite{Hu}, \cite{Pa}). A
similar formula holds for $C_{\mfl,j}$ and we conclude that
$C_{\mfl,j}$ and $C_{\mfk,j}$ are rational numbers. If there
exists a binormal Einstein metric on $M$, then $C_{\mfl,j}$ and
$C_{\mfk,j}$ are related by formula (\ref{caskjlj}) stated in this result. This implies that $\sqrt{2\al+1}$ is a rational
number.

$\Box$

\li

\bex  \textbf{Circle Bundles over Compact Irreducible Hermitian
Symmetric Spaces:} \eex
An application of Corollary \ref{bincor2} occurs when the base space is an irreducible symmetric
space. So let us consider a fibration $M \rightarrow
N$ where the fiber $F$ is isomorphic to the circle group and $N$ is an
isotropy irreducible symmetric space. Since $F$ is the circle
group, $\mfp$ lies in the center of $\mfk$. Hence, $K$ has
one-dimensional center, since for a compact irreducible symmetric
space the center of $K$ has at most dimension 1. Moreover, in this
case $N$ is a compact irreducible Hermitian symmetric space. In
particular, $L$ must coincide with the semisimple part of $K$. Compact irreducible Hermitian symmetric spaces $G/K$ are classified (see
e.g. \cite{He}). All the possible $G$, $K$ and $L$ are listed in Table \ref{sf2}, together with the coefficient $X$ of
the, unique, Einstein adapted metric on $G/L$, as in Corollary \ref{bincor2}.

$\Box$

\li

\btab[h!]\caption{Circle bundles over compact irreducible hermitian
symmetric spaces.}\label{sf2}
$$\bar{|lll|c|}\hline\xstrut  G & K & L & X \\ \hline

\xstrut SU(n) & S(U(p)\times U(n-p)) & SU(p)\times SU(n-p) &
\frac{p(n-p)+1}{2p(n-p)}\\

\xstrut SO(2n) &  U(n) &   SU(n) & \frac{n(n-1)+2}{2n(n-1)}\\

\xstrut SO(n) & SO(2)\times SO(n-2) & SO(n-2) & \frac{n-1}{n-2}\\

\xstrut Sp(n) & U(n) & SU(n) & \frac{n(n+1)+2}{2n(n+1)}\\

\xstrut E_6 & SO(10)\times U(1) & SO(10) & \frac{17}{32}\\

\xstrut E_7 & E_6\times U(1) & E_6 & \frac{14}{27}\\ \hline\ear$$
\etab

\section{Riemannian Fibrations with Einstein Fiber and Einstein Base}\label{einsteinFB}

In this section we investigate conditions for the existence of an
Einstein adapted metric $g_M$ on $M$ such that $g_F$ and $g_N$ are also Einstein. We prove Theorems \ref{bfein} and \ref{einrel} stated in Section \ref{intro}. We follow the notation and hypothesis introduced in Section \ref{intro}. In particular, for any homogeneous fibration $M=G/L\rightarrow G/K=N$, $G$ is a compact connected semisimple Lie group and
$L\varsubsetneq K\varsubsetneq G$ are connected closed non-trivial
subgroups of $G$.

\bprop \label{gnein} Let $g_M$ be an adapted metric on the homogeneous fibration $M=G/L\rightarrow G/K=N$ defined by the Euclidean product

$$g_{\mfm}=\left(\oplus_{a=1}^s\la_aB_{\mfp_a}\right)\oplus\left(\oplus_{k=1}^n\mu_kB_{\mfn_k}\right).$$

If $g_M$ and $g_N$ are both Einstein, then

$$\frac{\mu_j}{\mu_k}=\frac{r_j}{r_k}=\left(\frac{b^j}{b^k}\right)^{\frac{1}{2}},\, j,k=1,\ldots,n,$$

where $b^j$ is the eigenvalue of the operator
$\sum_{a=1}^s\la_aC_{\mfp_a}$ on $\mfn_j$, and
the $r_j$'s are determined by
$Ric^N=\oplus_{k=1}^nr_kB_{\mfn_k}$ as in Lemma
\ref{ricN}. Up to homothety, there exists at most one Einstein
metric $g_N$ on $N$ such that the corresponding adapted metric $g_M$ on $M$ is
Einstein.\eprop

\bproof The statement is trivial if $N$ is isotropy irreducible, so we suppose that $N$ is not irreducible, i.e., $n>1$.From Corollary \ref{cond1} we know that if
$g_M$ is Einstein, then there are constants $b^j$ such that

$$\sum_{a=1}^s\la_aC_{\mfp_a}\mid_{\mfn_j}=b^jId.$$

We recall from Lemma \ref{ricN} that
$Ric^N=\oplus_{k=1}^nr_kB_{\mfn_k}$. Hence, if $g_N$
is Einstein, then

\beq\label{rmus1}\frac{r_1}{\mu_1}=\ldots=\frac{r_n}{\mu_n},\textrm{  i.e., } \frac{\mu_j}{\mu_k}=\frac{r_j}{r_k}, \textrm{ for
}j,k=1,\ldots,n.\eeq

From Proposition \ref{riccin}, for $X\in\mfn_k$, the Ricci curvature
of $g_M$ is

$$Ric(X,X)=-\dfrac{1}{2\mu_k}\sum_{a=1}^s\la_aB(C_{\mfp_a}X,X)+r_kB(X,X)=\big(-\dfrac{b^k}{2\mu_k}+r_k\big)B(X,X).$$

If $g_M$ is Einstein, then from the expression above we obtain the
following Equations

\beq\label{rmus3}
-\dfrac{b^k}{2\mu_k^2}+\frac{r_k}{\mu_k}=-\dfrac{b^j}{2\mu_j^2}+\frac{r_j}{\mu_j}.\eeq

The identities (\ref{rmus1}) and (\ref{rmus3}) imply that

\beq\label{mubs}\dfrac{b^k}{\mu_k^2}=\dfrac{b^j}{\mu_j^2}\eeq

and consequently, by using (\ref{rmus1}),

$$\left(\frac{r_j}{r_k}\right)^2=\left(\frac{\mu_j}{\mu_k}\right)^2=\frac{b^j}{b^k}.$$

Since the eigenvalues $b^j$ are independent of the $\mu_j$'s, the ratios
(\ref{mubs}) imply that there is at most one possible choice for $g_N$, up to
scalar multiplication.

$\Box$\eproof

\bprop \label{gfein} Let $g_M$ be an adapted metric on the homogeneous fibration $M=G/L\rightarrow G/K=N$ defined by the Euclidean product
$$g_{\mfm}=\left(\oplus_{a=1}^s\la_aB_{\mfp_a}\right)\oplus\left(\oplus_{k=1}^n\mu_kB_{\mfn_k}\right).$$

If $g_M$ and $g_F$ are both Einstein, then

$$\frac{\la_a}{\la_b}=\frac{q_a}{q_b}= \sum_{j=1}^n\frac{C_{\mfn_j,b}}{\mu_j^2}\Big/\sum_{j=1}^n\frac{C_{\mfn_j,a}}{\mu_j^2} ,\textrm{ for
}a,b=1,\ldots,s,$$

where $c_{\mfn_j,a}$ is such that
$\Phi(C_{\mfn_j}\cdot,\cdot)\mid_{\mfp_a\times\mfp_a}=c_{\mfn_j,a}\Phi\mid_{\mfp_a\times\mfp_a}$,
for $a=1,\ldots,s$ and the $q_a$'s are determined by
$Ric^F=\oplus_{a=1}^sq_aB_{\mfp_a}$ as in Lemma
\ref{ricF}. Up to homothety, there exists at most one
Einstein metric $g_F$ on $F$ such that the corresponding adapted metric
$g_M$ on $M$ is Einstein.\eprop

\bproof The statement is trivial if $F=K/L$ is isotropy irreducible, so we suppose that $F$ is not irreducible, i.e., $s>1$. The proof is similar to that of Proposition \ref{gnein}, by
using Lemma \ref{ricF} and Proposition \ref{riccip}.

$\Box$\eproof

\textbf{\emph{Proof of Theorem \ref{bfein}:}} Using Proposition \ref{gnein}, we write
$\mu_j^2=\frac{b^j}{b^1}\mu_1^2$. The second formula follows
immediately from this and Proposition \ref{gfein}.

$\Box$

\li

\textbf{\emph{Proof of Theorem \ref{einrel}:}} If $g_M$, $g_N$ and $g_F$ are all Einstein,
from Proposition \ref{riccin} we get

$$-\frac{1}{2\mu_j}b^j+\mu_j E_N=\mu_j E ,$$

from which, if $E_N\neq E$, we deduce

\beq\label{eqmuj1}\mu_j^2=\frac{b^j}{2(E_N-E)}.\eeq

From Proposition \ref{riccip}, we have

$$\la_a E_F+\frac{\la_a^2}{4}\frac{C_{\mfn_j,a}}{\mu_j^2}=\la_a E$$

which implies

\beq \la_a\sum_{j=1}^n\frac{C_{\mfn_j,a}}{\mu_j^2}=4(E-E_F).\eeq

We obtain the required formula for $\la_a$ by replacing (\ref{eqmuj1}) in the
equation above.

$\Box$

\section{Invariant Einstein Metrics on Kowalski $n$-symmetric Spaces}\label{kow}

In this section we show that the $n$-symmetric spaces $M=\frac{G_0^n}{\triangle^nG_0}$ admit a non-standard Einstein adapted metric, for $n>4$. The main result is Theorem \ref{genknss} which is stated in section \ref{intro} and proved in this section.

Let $G_0$ be a compact connected simple Lie group. For any positive integer $m$, we denote $G_0^m=\underbrace{G_0\times\ldots\times G_0}_m$ and $\triangle ^mG_0$ is the diagonal subgroup in $G_0^m$.

Let $n,p,q$ be positive integers such that $p+q=n$ and $2\leq p\leq q\leq n-2$ and consider the following groups:

\beqar G = G_0^n\\ \nonumber
 K = \triangle^p G_0\times \triangle^q G_0 \varsubsetneq G\\ \nonumber
 L=\triangle ^nG_0\varsubsetneq K\eeqar

\li

We consider the fibration $M=G/L\rightarrow G/K=N$ with fiber $F=K/L$. If $\mfg_0$ denotes the Lie algebra of $G_0$, then the Lie algebras of $G$, $K$ and $L$ are $\mfg=\mfg_0^n$, $\mfk=\triangle^p \mfg_0\times \triangle^q \mfg_0$ and $\mfl=\triangle^n \mfg_0$, respectively. If $\Phi_0$ is the Killing form of $\mfg_0$, then the Killing form of $\mfg$ is simply $\Phi=\Phi_0+\ldots+\Phi_0$. Following the notation used in previous sections, let $\mfn$ be the orthogonal complement of $\mfk$ in $\mfg$ and
$\mfp$ be an orthogonal complement of $\mfl$ in $\mfk$, with respect
to $\Phi$. Then $\mfg=\mfl\oplus\mfm$, where $\mfm=\mfp\oplus\mfn$ and $\mfk=\mfl\oplus\mfp$.

\li

In the Lemma below we present a decomposition of the orthogonal complements $\mfp$ and $\mfn$. The tangent space to the fiber, $\mfp$, is $Ad\,L$-irreducible whereas $\mfn$ is $Ad\,K$-reducible. The decomposition below for $\mfn$ is among many others and the decomposition we will consider throughout. The proof is straightforward and can be found in \cite[\S 4.1]{Fa}.

\blem (i)
$\mfp=\{(\underbrace{qX,\ldots,qX}_p,\underbrace{-pX,\ldots-pX}_q):X\in\mfg_0
\}$ and $\mfp$ is $Ad\,L$-irreducible;

\li

(ii) $\mfn$ admits the decomposition $\mfn=\mfn_1\oplus\mfn_2$, where

\beqar\nonumber\mfn_1=\{(X_1,\ldots,X_p,0,\ldots,0):X_j\in\mfg_0,\sum_{j=1}^pX_j=0
\}\subset \mfg_0^p\times 0_q\\\nonumber
\mfn_2=\{(0,\ldots,0,X_1,\ldots,X_q):X_j\in\mfg_0,\sum_{j=1}^qX_j=0
\}\subset 0_p\times\mfg_0^q\eeqar

\elem

\li

Below we describe the Casimir operators of $\mfg$, $\mfk$, $\mfl$, $\mfp$, $\mfn_1$ and $\mfn_2$ and present the necessary eigenvalues to solve the Einstein equations for an Einstein adapted metric on $M$. The proofs of the following Lemmas can be found in \cite[\S 4.1]{Fa}.

\blem \label{ck1}(i) $C_{\mfg}=Id_{\mfg}$;

\li

(ii) $C_{\mfl}=\dfrac{1}{n}Id_{\mfg}$;

\li

(iii) $C_{\mfp}=\dfrac{q}{np}Id_{\mfg_0^p}\times
\dfrac{p}{nq}Id_{\mfg_0^q} $;

\li

(iv) $C_{\mfk}=\dfrac{1}{p}Id_{\mfg_0^p}\times
\dfrac{1}{q}Id_{\mfg_0^q} $;

\li

(v) $C_{\mfn_1}=\left(1-\dfrac{1}{p}\right)Id_{\mfg_0^p}\times
0_{\mfg_0^q}$ and  $C_{\mfn_2}=0_{\mfg_0^p}\times
\left(1-\dfrac{1}{q}\right)Id_{\mfg_0^q}$.
 \elem

We recall that $c_{\mfl,\mfp}$ is the eigenvalue of
$C_{\mfl}$ on $\mfp$, $c_{\mfk,i}$ is the eigenvalue of $C_{\mfk}$
on $\mfn_i$. The Casimir operator of $\mfp$ is scalar on $\mfn_i$, as we can see from Corollary \ref{ck2}, and $b^i$ denotes the eigenvalue of $C_{\mfp}$ on $\mfn_i$,
for $i=1,2$. Also, following a notation similar to that of (\ref{gas}) and (\ref{cnjs}), $c_{\mfn_i,\mfp}$ and $\ga$ are the constants
defined by

\beqar
\Phi(C_{\mfn_i}.,.)\mid_{\mfp\times\mfp}=c_{\mfn_i,\mfp}\Phi\mid_{\mfp\times\mfp},\,i=1,2\\
\Phi(C_{\mfk}.,.)\mid_{\mfp\times\mfp}=\ga \Phi\mid_{\mfp\times\mfp}.\eeqar

\blem \label{ck2}(i) $c_{\mfl,\mfp}=\dfrac{1}{n}$;

(ii) $C_{\mfp}$ is scalar on $\mfn_j$, $j=1,2$ and
$b^1=\dfrac{q}{np}$ and $b^2=\dfrac{p}{nq}$;

(iii) $c_{\mfk,1}=\dfrac{1}{p}$ and $c_{\mfk,2}=\dfrac{1}{q}$;

(iv) $\ga=\dfrac{q^2+p^2}{npq}$;

(v) $c_{\mfn_1,\mfp}=\dfrac{(p-1)q}{pn}$ and
$c_{\mfn_2,\mfp}=\dfrac{(q-1)p}{qn}$.

\elem

\li

We consider an adapted metric $g_M$ on $M$ defined by the Euclidean product

\beq\label{mdefknss}g_{\mfm}=\la B_{\mfp}\oplus\mu_1 B_{\mfn_1}\oplus\mu_2 B_{\mfn_2}.\eeq

 We observe that $\mfn_1$ and $\mfn_2$ are inequivalent
$Ad\,K$-modules, but they are not irreducible, for $n>4$. Hence, adapted
metrics on $M$ are not necessarily of the form (\ref{mdefknss}). However, throughout we shall focus only on adapted
metrics of this form.

We shall classify all the binormal Einstein metrics on $M$. Although the submodules $\mfn_1$ and $\mfn_2$ are not $Ad\,K$-irreducible for $n> 4$, the Casimir operators of $\mfl$, $\mfk$ and $\mfp$ are always scalar on  $\mfn_1$ and
on $\mfn_2$. Hence, it is enough to consider one irreducible submodule in
$\mfn_1$ and one irreducible submodule in $\mfn_2$, to compute the Ricci curvature of $g_M$. According to Theorem \ref{binormal1} (2.19), there is an
one-to-one correspondence, up to homothety, between binormal adapted
Einstein metrics on $M$ and positive solutions of the following
set of equations:

\beqar \label{eink1}\de_{12}^{\mfk}(1-X)=\de_{12}^{\mfl}\\
\label{eink2}(\ga+2c_{\mfl,\mfp})X^2-\left(1+2c_{\mfk,j}\right)X+(1-\ga+2b^j)=0,\,j=1,2\eeqar

Given a positive solution $X$, then Einstein binormal metrics are given by $g_{\mfm}=B_{\mfp}\oplus XB_{\mfn}$, up to homothety.

\bthm\label{binormalknss} Let $G_0$ be a compact  connected simple group and consider the fibration
$$M=\frac{G_0^n}{\triangle^nG_0}\rightarrow
\frac{G_0^p}{\triangle^pG_0}\times \frac{G_0^q}{\triangle^qG_0},$$

where $p+q=n$ and $2\leq p\leq q\leq n-2$.

If $p\neq q$ or $n=4$, the standard metric is the only Einstein binormal metric, up to homothety. For $n>4$ and $p=q$, there are on $M$ exactly two Einstein binormal metrics, up to homothety, which are the standard
metric and the metric induced by

$$g_{\mfm}=B_{\mfp}\oplus \frac{n}{4}B_{\mfn}.$$ \ethm

\bproof From Corollary \ref{ck2} we obtain that
$\de^{\mfk}_{12}=c_{\mfk_1}-c_{\mfk,2}=\frac{1}{p}-\frac{1}{q}$
whereas
$\de_{12}^{\mfl}=c_{\mfl_1}-c_{\mfl,2}=\frac{1}{n}-\frac{1}{n}=0$.
Hence, Equation (\ref{eink1}) implies that $X=1$ or $p=q$. So if
$p\neq q$, if there exists a binormal Einstein metric it must be the
standard metric. This we already know it is Einstein by \cite{Ro1}.
Therefore, if $p\neq q$, then there exists, up to homothety, exactly
one binormal Einstein metric on $G/L$ which is the standard metric.

By using Corollary \ref{ck2} Equation (\ref{eink2}) may be rewritten
as

\beqar\label{eink3}nX^2-q(p+2)X+pq+q-p=0,  \textrm{ for }j=1 \\
\label{eink4}nX^2-p(q+2)X+pq+p-q=0,  \textrm{ for }j=2 \eeqar

It is clear that $X=1$ is actually a solution of both equations, and we confirm that the standard metric is Einstein. Now suppose that $p=q$. As $n=p+q$, then $p=q=\frac{n}{2}$.
Therefore, (\ref{eink3}) and (\ref{eink4}) become equivalent to $4X^2-(n+4)X+n=0$. This polynomial has two positive roots, $1$ and $\frac{n}{4}$.
Therefore, for $p=q$ and $n>4$, there exist precisely two binormal
Einstein metrics.

$\Box$\eproof

Since the vertical isotropy space $\mfp$ is $Ad\,L$-irreducible, the restriction of any $G$-invariant metric on $M$ to the fiber $F$ is an Einstein metric. Below we show that the only Einstein adapted metric on $M$ which projects onto an Einstein metric on the base space $N$ is the Einstein binormal metric given in Theorem \ref{binormalknss}.

\bthm \label{form2} Let $G_0$ be a compact  connected simple group and consider the fibration
$$M=\frac{G_0^n}{\triangle^nG_0}\rightarrow
\frac{G_0^p}{\triangle^pG_0}\times \frac{G_0^q}{\triangle^qG_0},$$

where $p+q=n$ and $2\leq p\leq q\leq n-2$. Let $g_M$ be an Einstein adapted metric on $M$
defined by $g_{\mfm}=\la B_{\mfp}\oplus\mu_1 B_{\mfn_1}\oplus\mu_2 B_{\mfn_2}$ as in (\ref{mdefknss}). The projection $g_N$ onto the
base space is also Einstein if and only if $p=q$. In this case, $g_M$ is a binormal metric. \ethm

\bproof By Proposition \ref{gnein} we know that if $g_M$ and $g_N$ are
Einstein then we must have the relation

\beq\label{rbsk}\frac{r_1}{r_2}=\left(\frac{b^1}{b^2}\right)^{\frac{1}{2}}.\eeq

From Lemma \ref{ck2} (ii) we obtain
$\left(\frac{b^1}{b^2}\right)^{\frac{1}{2}}=\frac{q}{p}$. Since
$[\mfn_1,\mfn_2]=0$ and $C_{\mfl}$ is scalar on $\mfn_i$, it follows from the definition of the $r_i$'s in Proposition \ref{riccin} and from Corollary \ref{ric1cor2} that $r_i=\frac{1}{2}\left(\frac{1}{2}+c_{\mfk,i}\right)$. Hence, $\frac{r_1}{r_2}=\frac{(p+2)q}{(q+2)p}$, by using Lemma
\ref{ck2} (iii). Therefore, (\ref{rbsk}) is possible if and only if
$p=q$. Also from the proof of Proposition \ref{gnein}, if $g_N$ and
$g_M$ are Einstein, then
$\frac{\mu_1}{\mu_2}=\left(\frac{b^1}{b^2}\right)^{\frac{1}{2}}=\frac{p}{q}=1$
and $g_M$ is binormal. Conversely, if $g_M$ is binormal and $p=q=n/2$, then by the above we also get

$$\frac{r_1}{\mu_1}=\frac{r_2}{\mu_2}$$ and $g_N$ is Einstein.

$\Box$\eproof

The Einstein equations in general for arbitrary $p$ and $q$ are extremely complicated. However with the help of Maple it is still possible to solve the problem in general. Next we shall classify all the Einstein adapted metrics on $M$ of the form (\ref{mdefknss}). The proof of the following Lemma follows from Proposition \ref{riccinf} and Corollary \ref{ric1cor2} by using the eigenvalues given in Corollary \ref{ck2}. This can be found in detail in \cite[\S 4.1]{Fa}.

\blem \label{knsseeq}Consider the fibration
$$M=\frac{G_0^n}{\triangle^nG_0}\rightarrow
\frac{G_0^p}{\triangle^pG_0}\times \frac{G_0^q}{\triangle^qG_0}=N,$$

where $p+q=n$ and $2\leq p\leq q\leq n-2$. There is a one-to-one
correspondence between Einstein adapted metrics on $M$ defined by $g_{\mfm}=\la B_{\mfp}\oplus\mu_1 B_{\mfn_1}\oplus\mu_2 B_{\mfn_2}$ as in (\ref{mdefknss}), up to homothety, and positive
solutions of the following system of Equations:

\beqar -2q^2X_1^2+nq(p+2)X_1+2p^2X_2^2-np(q+2)X_2=0 \\
n^2+q^2(p+1)X_1^2+p^2(q-1)X_2^2-nq(p+2)X_1=0 \eeqar

To a positive solution $(X_1,X_2)$ corresponds an Einstein adapted metric
defined by $g_{\mfm}=B_{\mfp}\oplus\frac{1}{X_1} B_{\mfn_1}\oplus\frac{1}{X_2} B_{\mfn_2}$.\elem

\li

\textbf{\emph{Proof of Theorem \ref{genknss}:}} By using Maple we obtain that the solutions of the system given in Lemma \ref{knsseeq} are $X_1=X_2=1$ and

\beq \label{solgen2}X_1=\alpha,\,X_2=\left(\frac{-q^2(p+1)\alpha^2+nq(p+2)\alpha-n^2}{p^2(q-1)}\right)^{\frac{1}{2}},\eeq

where $\alpha$ is a root of the polynomial

$$t(Z)=4q^2Z^3-4q(n+pq+2)Z^2+n(q(q+2)(p+1)+n+8)Z-(q+3)n^2.$$

The solution $X_1=X_2=1$  corresponds to a standard metric and, once more, we confirm that $M$ is an Einstein standard manifold. We investigate the existence of other metrics. From the expression for $X_2$ in (\ref{solgen2}) we conclude that,

$$X_2\in\reals\textrm{ if and only if }\al\in\left(\frac{n}{q(p+1)},\frac{n}{q}\right).$$

For this we compute the roots of the polynomial $-q^2(p+1)\alpha^2+nq(p+2)\alpha-n^2$ in (\ref{solgen2}). Simple calculations show that

\beqar t\left(\frac{n}{q}\right)=\frac{p(q-1)^2n^2}{q}>0\nonumber\\
t\left(\frac{n}{q(p+1)}\right)=-\frac{p(p+3)^2(q-1)n^2}{q(p+1)^3}<0\nonumber
\eeqar

and thus $t$ has at least one (positive) root in the interval $\left(\frac{n}{q(p+1)},\frac{n}{q}\right)$. To this root corresponds an Einstein adapted metric on $M$. Furthermore, we show that this root is unique and distinct from 1. From this we will conclude that there exists a non-standard Einstein adapted metric on $M$. Simple calculations show that the zeros of $\frac{dt}{dZ}$ are

$$\frac{n+pq+2}{3q}\pm\frac{\sqrt{\de}}{6q},$$

where $$\de=(q+1)^2p^2-(q-1)(3q^2+4q-8)p-(q-1)(3q^2+8q+16).$$

We show that $\de<0$. For $p=q$, $\de=-2q^4-2q^3+8q^2-16q+16<0$, for every $q\geq 2$. So we suppose that $p<q$. In this case, since $p^2\leq (q-1)p$, we have

$$\bar{rl}\de\leq & (q-1)\big(-(2p+3)q^2-(2p+8)q+(9p-16)\big)\\

< & (q-1)(-2p^3-5p^2+p-16)\\

< & 0,\ear$$

for every $p\geq 2$. With this we conclude that $\frac{dt}{dZ}$ has no real zeros and thus the root of $t$ found above is the unique real root of $t$. Moreover, we must guarantee that this root does not yield the solution $X_1=X_2=1$. If $X_1=X_2=1$, then $\al=1$ is a root of $t$. This may be possible since $1\in\left(\frac{n}{q(p+1)},\frac{n}{q}\right)$. Since

$$t(1)=p(q+2)(q-1)(n-4),$$

 $\al=1$ is a root of $t$ if and only if $n=4$.
By using (\ref{solgen2}) we get that if $X_1=1$ when $n=4$, then
$X_2=1$ as well. Since non-standard Einstein adapted metrics are
given by pairs of the form (\ref{solgen2}), with $\al\neq 1$, we
conclude that there exists a unique non-standard Einstein adapted
metric of the form (\ref{mdefknss}) if and only if $n>4$; in
the case $n=4$, the standard metric is the unique Einstein adapted
metric of the form (\ref{mdefknss}). Finally, we observe
that, if $n=4$, the subspaces $\mfn^1$ and $\mfn^2$ are
irreducible $Ad\,L$-submodules. Hence, any adapted metric on $M$
is of the form (\ref{mdefknss}). Therefore, we conclude
that, for $n=4$, there exists a unique Einstein adapted metric on
$M$ which is the standard metric.

Since there is a unique non-standard Einstein adapted metric on
$M$, it follows from Theorem \ref{binormalknss} that this metric
is binormal if and only if $p=q$.

$\Box$


\begin{thebibliography}{99}

\bibitem{ADF}D. Alekseevsky, I. Dotti, C. Ferraris (1996) \textit{Homogeneous Ricci Positive $5$-Manifolds}, Pacific Journal of Mathematics, 175 (1), 1-12.

\bibitem{AP}D. Alekseevsky, A.M. Perelomov (1987) \textit{Invariant K\"{a}hler-Einstein Metrics on Compact Homogeneous Spaces}, Moscow Correspondence State Pedagogical Institute, Institute of Theoretical and Experimental Physics; translated from Funktional'nyi Analiz i Ego Prilozheniya, 20 (3), 1-16.

\bibitem{Fa}F. Ara\'{u}jo (2008) \textit{Einstein Homogeneous Riemannian Fibrations}, PhD Thesis of The University of Edinburgh, arXiv:0905.3143v1 [math.DG].

\bibitem{BB}L. B\'{e}rard-Bergery (1975) \textit{Sur Certaines
Fibrations d'Espaces Homogene Riemanniens}, Compositio Math. 30,
43-61.

\bibitem{Be}A. Besse (1987) \textit{Einstein Manifolds},
Springer-Verlag.

\bibitem{BK}C. B\"{o}hm, M.M. Kerr (2005) \textit{Low-dimensional Homogeneous Einstein Metrics}, Trans. Amer. Math. Soc. 358 (4), 1455-1468.

\bibitem{BG}C.P. Boyer, K. Galicki (2000) \textit{On Sasakian-Einstein Geometry}, Internat.J.Math. 11, 873-909.

\bibitem{CRW}L. Castellani, L.J. Romans, N.P. Warner (1984) \textit{A Classification of Compactifying Solutions for d=11 Supergravity}, Nuclear Physics 241 B, 429-262.

\bibitem{DK}W. Dickinson, M.M. Kerr (2008) \textit{The Geometry of Compact Homogeneous Spaces with Two Isotropy Summands}, Ann Glob Anal Geom. 34, 329-350.

\bibitem{He}S. Helgason (1978) \textit{Differential Geometry, Lie Groups and
Symmetric Spaces}, Ac.Press, Mathematics 80.

\bibitem{Hu}J.E Humphreys (1972) \textit{Introduction to Lie Algebras and Representation
Theory}, Springer-Verlag.

\bibitem{Je1}G.R. Jensen (1969) \textit{Homogeneous Einstein Spaces of Dimension 4}, J.Diff.Geometry, 3, 309-349.

\bibitem{Je2}G.R. Jensen (1973) \textit{Einstein Metrics on Principal Fiber Bundles}, J.Diff.Geometry, 8, 599-614.

\bibitem{Ke1}M.M. Kerr (1996) \textit{Some New Homogeneous Einstein
Metrics on Symmetric Spaces}, Trans. Amer. Math. Soc. 348 (1),
153-11.

\bibitem{NK}S. Kobayashi, K. Nomizu (1969) \textit{Foundations of Differential Geometry I, II}, Tracts in Mathematics 15.

\bibitem{Ko}O. Kowalski (1980) \textit{Generalized Symmetric Spaces}, Springer-Verlag.

\bibitem{Ni}Y.G. Nikonorov (2004) \textit{Compact Homogeneous Einstein 7-manifolds}, Geom. Dedicata 109, 7-30.

\bibitem{NR2}Y.G. Nikonorov, D. E. Rodionov (1999) \textit{Compact 6-dimensional Homogeneous Einstein Manifolds}, Dokl.Math. RAR 336, 599-601.

\bibitem{No}K.Nomizu (1954) \textit{Affine Connections on Homogeneous Spaces}, Amer.J.Mathematics 76 (1), 33-65.

\bibitem{Pa}D. Panyushev (2001) \textit{Isotropy Representations, Eingenvalues of a Casimir Element and Commutative Lie subalgebras},
J.London Math. Soc. 64 (2), 61-80.

\bibitem{Ro1}E.D. Rodionov (1995) \textit{Homogeneous Riemannian Manifolds
With Einstein Metrics}, Proc.conf., Aug 28-Sep 1, Brno Czech Republic, Masaryk University, 81-91.

\bibitem{Ro4}E.D. Rodionov (1992) \textit{simply-connected Compact Standard Homogeneous Einstein Manifolds}, Siberian Mathematical
Journal, 33 (4), 104-119.

\bibitem{Sa}A.A. Sagle (1970) \textit{Some Homogeneous Einstein
Manifolds},  Nagoya Math. J. 39, 81-106.

\bibitem{Ti}G. Tian (1997) \textit{K\"{a}hler-Einstein Metrics with Positive Scalar Curvature}, Invent. Math. 130, 1-37.

\bibitem{WZ}M.Y. Wang, W.Ziller (1985) \textit{On Normal Homogeneous Einstein
Manifolds}, Annales Scientifiques de l'E.N.S., 18 (4), 563-633.

\bibitem{WZ2}M.Y. Wang, W. Ziller (1986) \textit{Existence
and Non-existence of Homogeneous Einstein Metrics}, Inventiones
Math. 84, 177-194.

\bibitem{WZ3}M.Y. Wang, W. Ziller (1990) \textit{Einstein Metrics on Principal Torus Bundles}, J.Diff.Geom. 31, 215-248.

\bibitem{Wo}J.A. Wolf (1968) \textit{The Geometry and Structure of
Isotropy Irreducible Homogeneous Spaces}, Acta Math. 120, 59-148;
Correction: 152, 141-142 (1984).

\bibitem{Yau}S.T. Yau (1978) \textit{On the Ricci Curvature of a Compact K\"{a}hler Manifold and the Complex Monge-Amp\`{e}re Equation I}, Comm. Pure Appl. Math. 31, 339-411.

\bibitem{Zi}W. Ziller (1984) \textit{Homogeneous Einstein Metrics}, Global
Riemannian Geometry, TJ.Willmore - N.J.Hitchin Eds, John Wiley,
126-135.

\bibitem{Zi2}W. Ziller (1982) \textit{Homogeneous Einstein Metrics on Spheres and Projective Spaces}, Math. Ann. 259, 351-358.

\end{thebibliography}
\end{document}